\documentclass[12pt]{article}

\usepackage[cp1251]{inputenc}
\usepackage[english]{babel}

\usepackage{amsfonts,epsfig,epstopdf}
\usepackage{amssymb,amsmath,enumerate,float}
\makeatother

\addto\captionsenglish{}

\renewcommand{\i}{i}
\newcommand{\ptl}{\partial}

\newcommand{\be}{\begin{equation}}
\newcommand{\ee}{\end{equation}}
\newcommand{\beq}{\begin{equation}}
\newcommand{\eeq}{\end{equation}}

\newcommand{\eps}{\varepsilon}

\title{Transient phenomena in a three-layer waveguide and the analytical 
structure of the dispersion diagram}
\author{A.~V.~Shanin, K.~S.~Knyazeva}

\begin{document}
\maketitle

\begin{abstract}
Excitation of waves in a three-layer acoustic wavegide is studied.
The wave field is presented as a sum of integrals. The summation is held over all 
waveguide modes. The integration is performed over the temporal frequency axis. 
The dispersion diagram of the waveguide is analytically continued,  
and the integral 
is transformed by deformation of the integration contour into the domain of complex frequencies.   
As the result, the expression for the fast components of the signal (i.~e.\ for the transient 
fields) is simplified. 

The structure of the Riemann surface of the dispersion diagram of the waveguide is studied.
For this, a family of auxiliary problems indexed by the parameters describing the 
links between layers is introduced. The family depends on the linking parameters analytically, 
and the limiting case of weak links can be solved analytically.
\end{abstract}

{\bf Keywords:} Waveguide, dispersion diagram, precursor, forerunner, firs-arriving signal, 
analytical continuation, orthogonality relations

\section{Introduction}

Let a layered waveguide be excited by a short pulse. A complicated wave process 
starts. In a very short time after the pulse one can expect a ray process 
near the source. After a very long time there should be a ``far-field'' stage. A comprehensive understanding of this stage is given by the dispersion diagram of the waveguide.
The wave components propagate with group velocities provided by the 
dispersion diagram. Between the ray phase and the far-field phase there can exist 
some transient waves, whose velocities may differ from the group velocities. If such 
transient waves are faster than the components of the far-field, they are called precursors, forerunners, 
of first--arriving signals (FAS). 

There are well-known types of precursors, namely the Sommerfeld's precursors and the Brillouin's precursors \cite{Brillouin,Akhmanov}. The Sommerfeld's precursors can be observed when the signal has high-frequency components traveling with a high velocity. The 
Brillouin's precursors can be observed when the spectrum of the excitation is not smooth.   

Here we study the case when the excitation has no high-frequency components, and it is 
smooth. Thus, the Sommerfeld's and Brillouin's precursors are not considered here. Instead, 
we study precursors of the type of leaky waves. Some general ideas for describing the transient waves of this type can be found in \cite{Shanin},
which is a development of \cite{Miklowitz}.   

The previous study of the leaky-wave type precursors demonstrated that an efficient description of them can be obtained by an analytic continuation of the dispersion diagram. Namely, the representation of the field is a series--integral expression. The series is
taken over the modes of the waveguides, while the integral is taken over the frequency axis ~$\omega$. The key idea is to deform the integration contour. While the contour is deformed, 
it crosses some branch points of the dispersion diagram. After such crossings the structure of the dispersion diagram on the contour becomes simpler. 

Thus, the structure of the Riemann surface of the dispersion diagram (the multivalued function $k(\omega)$) is of considerable importance for the description of transient properties in waveguides. Generally, the structure of the Riemann surface is not known. 
In the current paper we describe the structure of the Riemann surface for a non-trivial case 
of a three-layer waveguide. The tool for studying this Riemann surface is a gradual 
switching on the links between the layers. If the links are weak, the structure of the Riemann surface is found from an asymptotic consideration. While the link becomes stronger, 
the surface changes homotopically, in particular, the branch points travel along some 
continuous trajectories in the $\omega$ plane. This evolution can be tracked numerically. Thus, 
we bring some order into the structure of the Riemann surface of the dispersion diagram and
provide a numerical technique for finding the positions of its branch points.     

The structure of the paper is as follows. In Section~\ref{sec_formulation} the problem for the 
waveguide is formulated. A dispersion equation for the waveguide is 
built. 

In Section~\ref{sec_miklowitz} the main ideas of the method proposed in 
\cite{Miklowitz,Randles,Shanin} are listed briefly. Namely, it is explained why it is important to study the 
analytical continuation of the dispersion diagram.  
A dispersion diagram is built numerically for some physical realization of the waveguide.
Besides the diagram for real $\omega$, some analytical continuations are demonstrated. 
It is shown that the real diagram has a terraced structure, i.~e.\ it has pseudo-crossings 
of branches, while the analytical continuations have crossings. The types of waveguide modes 
related to different branches of the diagram (real or continued) are discussed. 
A numerical demonstration of the Miklowitz--Randles method 
is presented. It is shown that fast components of the signals (the precursors) can be described by 
a small amount of terms in the series-integral representation after the contour is deformed. 

in Section~\ref{sec_branchpoints} a method to study the Riemann surface of the dispersion diagram 
is proposed. The initial problem is embedded into an analytical family of 
auxiliary problems depending on the linking parameters. The Riemann surface of the problem corresponding to small linking parameters 
can be easily constructed. A numerical method of finding the trajectories of the branch points as 
the linking parameters grow from small values to infinity is proposed and validated. 
In Appendix~A a set of bilinear relations for the waveguide modes is constructed. 
The most important of them is the orthogonality relation for the mode corresponding 
to a branch point of the dispersion diagram. In Appendix~B the branch points positions
are found for small linking parameters using the perturbation method.

\section{Problem formulation}
\label{sec_formulation}

\subsection{Equations and boundary conditions}

\begin{figure}[ht]
\centerline{\epsfig{file=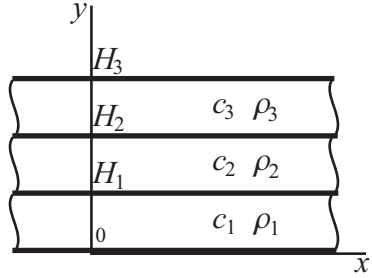}}
\caption{Geometry of the waveguide}
\label{fig01}
\end{figure}

Consider a planar waveguide consisting of three liquid layers. In the $(x,y)$--plane the 
layers~1, 2, 3 occupy the strips $0<y<H_1$, $H_1 < y < H_2$, $H_2 < y < H_3$, respectively 
(see Fig.~\ref{fig01}). 
The densities of the media are equal to $\rho_1$, $\rho_2$, $\rho_3$. The sound velocities 
are equal to $c_1$, $c_2$, $c_3$.     
The thicknesses of the layers are 
\[
h_1 = H_1, \qquad h_2 = H_2 - H_1, \qquad h_3 = H_3-H_2 . 
\]
It is convenient to define the functions
\[
c(y) = \left\{ \begin{array}{ll}
c_1, & 0 < y < H_1 \\
c_2, & H_1 < y < H_2 \\
c_3, & H_2 < y < H_3
\end{array} \right.
\qquad
\rho (y) = \left\{ \begin{array}{ll}
\rho_1, & 0 < y < H_1 \\
\rho_2, & H_1 < y < H_2 \\
\rho_3, & H_2 < y < H_3
\end{array} \right.
\]

Acoustic field in the layers is described by the acoustic potential $u (x,y,t)$
linked with the pressure and the particle velocity by the relations 
\begin{equation}
p = \rho(y) \frac{\ptl u}{\ptl t},
\qquad 
{\rm v}_j = - \nabla u_j.
\label{eq0101}
\end{equation} 
The potential obeys the wave equation: 
\begin{equation}
\Delta u - \frac{1}{c^2} \frac{\ptl^2 u}{\ptl t^2} = 0 , 
\label{eq0102}
\end{equation}
which is valid for $y \ne H_1, H_2$.

If an excitation of the waveguide is considered, an inhomogeneous 
Helmholtz equation is introduced:  
\begin{equation}
\Delta u - \frac{1}{c^2} \frac{\ptl^2 u}{\ptl t^2}= \delta(x) \delta(y-y_0) f(t), 
\label{eq0104}
\end{equation}
where $\delta$ is the Dirac delta-function, $f(t)$ is the excitation profile.  

Assume that at the external boundaries of the waveguide are acoustically hard (Neumann):
\begin{equation}
\left. \frac{\ptl u }{\ptl y} \right|_{y = 0}= 0, 
\qquad 
\left. \frac{\ptl u }{\ptl y} \right|_{y = H_3}= 0. 
\label{eq0103}
\end{equation}

On the internal boundaries the normal velocity and the pressure should be continuous. 
Thus, the boundary conditions are as follows: 
\begin{equation}
\left. \frac{\ptl u }{\ptl y} \right|_{y = H_1-0} = 
\left. \frac{\ptl u }{\ptl y} \right|_{y = H_1+0} ,
\qquad 
\left. \frac{\ptl u }{\ptl y} \right|_{y = H_2-0} = 
\left. \frac{\ptl u }{\ptl y} \right|_{y = H_2+0} ,
\label{eq0105}
\end{equation}
\begin{equation}
\left. \rho_1 \frac{\ptl u }{\ptl t} \right|_{y = H_1-0} = 
\left. \rho_2 \frac{\ptl u }{\ptl t} \right|_{y = H_1+0} ,
\qquad 
\left. \rho_2 \frac{\ptl u }{\ptl t} \right|_{y = H_2-0} = 
\left. \rho_3 \frac{\ptl u }{\ptl t} \right|_{y = H_2+0} .
\label{eq0106}
\end{equation}


\subsection{Dispersion diagram}
\label{sec_diagram}

The formulation above is related to a time-dependent process in the waveguide. To study 
free wave propagation, however, it is necessary to consider the process that is stationary 
(harmonic) in time and is exponential with respect to $x$-coordinate. Let the  
temporal circular frequency be equal to $\omega$ and the $x$-wavenumber to $k$, i.~e. consider 
the representation 
\begin{equation}
u(x,y,t) = \exp \{ i k x - i \omega t \} u (y).
\label{eq0107}
\end{equation}
The wave equation (\ref{eq0108})  becomes a set of ordinary differential equations 
\begin{equation}
\frac{d^2 u}{d y^2} + \alpha_j^2 u = 0 ,
\qquad
\alpha_j = 
\alpha_j (\omega, k)=
\sqrt{\frac{\omega^2}{c_j^2} - k^2}, 
\label{eq0108}
\end{equation}
which are valid in the strips $0<y<H_1$, $H_1<y<H_2$, $H_2<y<H_3$ for $j = 1,2,3$, respectively.

The signs of the radicals are not important (the final equations do not depend on these signs), 
but for definiteness we use the standard convention  
taking positive real or positive imaginary values for real $\omega$ and~$k$.
Boundary conditions (\ref{eq0103}), (\ref{eq0105}), and (\ref{eq0106})
become as follows: 
\begin{equation}
u'(0)= 0, 
\qquad 
u'(H_3)= 0, 
\label{eq0109}
\end{equation}
\begin{equation}
u' (H_1-0) = 
u' (H_1+0) ,
\qquad 
u' (H_2-0) = 
u' (H_2+0), 
\label{eq0110}
\end{equation}
\begin{equation}
\rho_1 u(H_1-0) = 
\rho_2 u(H_1+0) ,
\qquad 
\rho_2 u(H_2-0) = 
\rho_3 u(H_2+0) 
.
\label{eq0111}
\end{equation}
Here prime denotes the $y$-derivative.

Consider the solution in the following form: 
\begin{equation}
u(y) = \left\{ \begin{array}{ll}
A \cos(\alpha_1 y),                      & 0<y<H_1 \\  
B \sin(\alpha_2 y) + C \cos(\alpha_2 y), & H_1<y<H_2 \\
F \cos(\alpha_3 (y-H_3)),                & H_2<y<H_3
\end{array} \right.
\label{eq0112}
\end{equation}
Such a choice makes conditions (\ref{eq0109}) be obeyed automatically. Then, from conditions 
(\ref{eq0110}) and (\ref{eq0111}) we obtain that a nontrivial solution of the problem exists if and only if the following dispersion relation is valid:
\begin{equation}
D(\omega,k) = 0,
\label{eq0113}
\end{equation}
\[
D(\omega,k) = {\rm det}({\bf D}(\omega, k)),
\]
\begin{equation}
{\bf D} = \left( \begin{array}{cccc}
-\rho_1 \cos (\alpha_1 H_1)  & \rho_2 \sin(\alpha_2 H_1) & \rho_2 \cos (\alpha_2 H_1) & 0  \\
\alpha_1 \sin (\alpha_1 H_1) & \alpha_2 \cos(\alpha_2 H_1) & -\alpha_2 \sin (\alpha_2 H_1) & 0  \\
0 & \rho_2 \sin (\alpha_2 H_2) & \rho_2 \cos(\alpha_2 H_2) & - \rho_3 \cos (\alpha_3 h_3) \\
0 & \alpha_2 \cos (\alpha_2 H_2) & 
-\alpha_2 \sin(\alpha_2 H_2) &  -\alpha_3 \sin (\alpha_3 h_3) 
\end{array} \right)
\label{eq0114}
\end{equation}
It is quite clear how (\ref{eq0113}) can be obtained. The representations (\ref{eq0112})
should be substituted into (\ref{eq0110}), (\ref{eq0111}). This forms four equations for four 
coefficients $A$, $B$, $C$, $F$. Relation (\ref{eq0113}) provides the existence of a non-trivial 
solution for this system. 

Instead of variables $(\omega, k)$ we find it convenient sometimes to use the variables 
\[
W = \omega^2, 
\qquad
K = k^2.
\]
The same determinant depending on variables $W,K$ will be denoted by 
\[
\tilde D (W,K) = D(\sqrt{W}, \sqrt{K}).
\]


\section{Miklowitz--Randles method}
\label{sec_miklowitz}

\subsection{Main ideas of Miklowitz--Randles method}

Let $k_n (\omega)$ be a set of all roots of (\ref{eq0113}) solved with respect to 
the variable~$k$. For each real $\omega$ there is an infinite set of roots $k_n$, 
however only a finite number of them are real. All other roots are imaginary. The real 
roots correspond to propagating wave modes, while the imaginary roots correspond to 
evanescent waves. 

Let the waveguide be excited by a source having the profile $f(t)$ 
(see (\ref{eq0104})). 
Let the observation point have coordinates $x = L$, $y = H_3$.    
The field $u$ at the observation point can be written as a modal expansion
of the ``complex'' signal:
\[
u(L,H_3,t) = {\rm Re} [\tilde u], 
\]
\begin{equation}
\tilde u (L, H_3, t) = \sum_n \int \limits_{0}^{\infty}
\hat f(\omega) Y(\omega, k_n) \exp \{ i k_n (\omega) L - i \omega t \} 
d\omega , 
\label{eq0201}
\end{equation}
where the summation is taken over all roots of the dispersion equation, 
$\hat f(\omega)$ is the Fourier 
transform of~$f$:
\begin{equation}
\hat f(\omega) = \frac{1}{2\pi} \int \limits_{-\infty}^{\infty}
f(t) \exp \{ i \omega t \} \, d t, 
\label{eq0202}
\end{equation} 
$Y (\omega, k_n)$ is a trigonometric function not containing $L$ (see below). 
The derivation of (\ref{eq0201}) is rather straightforward. This derivation can be found, for example, 
in \cite{Shanin} for the case of a two--layer waveguide. The three--layer case is quite similar. 

The analysis of the transient waves in a layered waveguide is based on the following statements. 
The listed ideas have been proposed in \cite{Shanin}, which in its turn is deeply based on 
\cite{Miklowitz}. 

\begin{figure}[ht]
\centerline{\epsfig{file=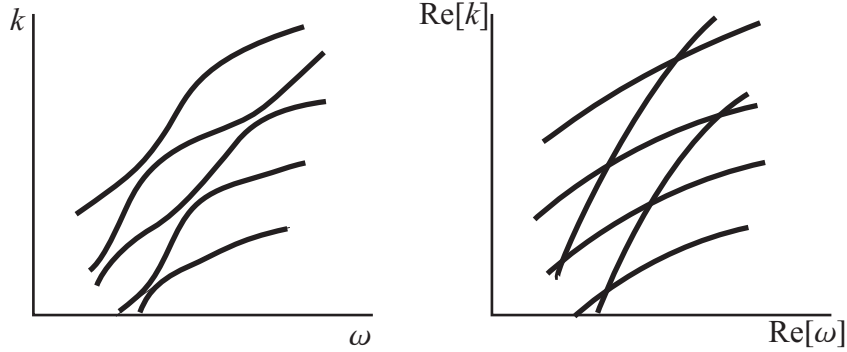}}
\caption{Terrace-like structures of the dispersion diagram (left); simplification of the 
diagram after analytical continuation (right)}
\label{fig02}
\end{figure}

{\bf 1.} The branches of the dispersion diagram $k_n (\omega)$ form a terrace-like structure,
see Fig.~\ref{fig02}, left.
(the term has been introduced in \cite{Mindlin} for an elastic waveguide). When there is 
a large amount of propagating modes, when there is a wide-band excitation pulse, and when 
$L$ is relatively small, 
the sum-integral
expression (\ref{eq0201}) is too complicated and cannot be asymptotically evaluated directly. So the expression needs some further transformations. Only after such transformation one can extract components that can be attributed as precursors.    

{\bf 2.} 
The functions $k_n (\omega)$ should be considered as a single {\em multivalued\/} function
$k(\omega)$
of the {\em complex\/} variable~$\omega$. The number of sheets of this function is infinite. The function has its Riemann surface in the usual sense. A direct analogy of this step is the well-known
mathematical approach to algebraic equations. The Galois theory studies all roots of an algebraic equation together. The branch points of the Riemann surface and the permutations of the roots happening at the branch points play an important role in  the theory. 

Here we have a transcendental equation (\ref{eq0113}) with respect to $k$, so we cannot expect any analytic representation for~$k(\omega)$.

{\bf 3.}
In the general case (random values of sound velocities) all branch points are of order~2, i.~e.\ connect two sheets. This fact can be explained as follows. 
Since $k(\omega)$ is set by an implicit relation $D(\omega , k(\omega)) =0$, the 
derivative of $k(\omega)$ is set by the relation 
\begin{equation}
\frac{d k}{d \omega} = \frac{\ptl D / \ptl \omega}{\ptl D / \ptl k}
\label{eq0203}
\end{equation}  
A point $(\omega_*, k_* = k(\omega_*))$ 
can be a singular point of the function $k(\omega)$ only if  
the denominator is zero, i.~e.\ if two
equations are fulfilled simultaneously: 
\begin{equation}
D(\omega_*, k_*) = 0, 
\qquad 
\ptl_k D(\omega_*, k_*) = 0, 
\label{eq0204}
\end{equation}
where $\ptl_k$ denotes a partial derivative with respect to the second argument. 
The system (\ref{eq0204}) is explained and studied in \cite{Randles} for the case of an elastic 
planar waveguide. 
If $(\omega_* , k_*)$ is such a point, it would be natural to suppose that 
near this point
\begin{equation}
D(\omega , k) = a (\omega- \omega_*) + b (k - k_*)^2 + c (\omega- \omega_*)(k - k_*) +
 \dots
\label{eq0205}
\end{equation} 
for some constants $a$, $b$, and $c$. The dots denote the smaller terms of the Taylor series. 
The expansion (\ref{eq0205}), when $D=0$ is solved with respect to $k$, leads to a square root 
behavior of $k(\omega)$, i.~e.\ to a branch point of the second order:
\[
k = k_* - \frac{1}{2b} \left(\pm 2i \sqrt{a b} (\omega- \omega_*)^{1/2} + 
c (\omega- \omega_*) + O((\omega- \omega_*)^3/2) \right) 
\]
A simple consideration shows that function $k(\omega)$ cannot have poles.  

The cut-off frequencies for the modes are also branch points of the dispersion diagram. 
They belong to the real axis of $\omega$ and connect the sheets corresponding to similar modes 
traveling forward and backward. These branch points do not play any important role 
in our consideration. We are interested mainly in the branch points located in the
upper half-plane. These branch points  ``connect'' different modes traveling in the positive direction.   

{\bf 4.}
The contour in (\ref{eq0201}) can be deformed {\em without taking into account the position of the 
branch points\/}, i.~e.\ passing the branch points if necessary. 
Such a deformation does not change the value of the sum of the integrals. 
This statement may sound surprising, but it can be easily explained. The integral 
in (\ref{eq0201}) is taken over all sheets of a Riemann surface of~$k(\omega)$. Assume that 
the contour is deformed continuously, and 
a branch point of order two is hit. According to the Cauchy's theorem, one should keep an additional 
loop encircling the branch point. However, such a loop appears on two sheets, which are connected at the branch point, and these two additional loops compensate each other. A detailed explanation 
can be found in \cite{Miklowitz} and in \cite{Shanin} (see Fig.~10 in \cite{Shanin} describing the rebuilding of branches of the dispersion diagram).

\begin{figure}[ht]
\centerline{\epsfig{file=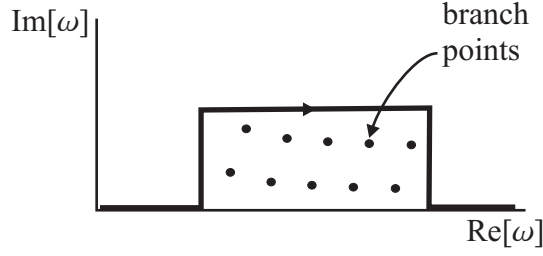}}
\caption{Local contour deformation for (\ref{eq0201})}
\label{fig03}
\end{figure}

{\bf 5.}
Under some conditions the contour of integration in (\ref{eq0201}) can be deformed as 
it is shown in Fig.~\ref{fig03}.
A part of the contour is shifted into the domain of positive imaginary $\omega$.
 The band within which the deformation is performed corresponds 
to the frequency band of the exciting pulse. The aim of the deformation is to pass 
the branch points connecting the sheets of the Riemann surface bearing the propagating modes. 
As the result, on the deformed contour the dispersion diagram becomes simplified 
(see Fig.~\ref{fig02}, right). 
Instead of a terrace-like structure, one gets two (or more) regular families of branches. The 
expression (\ref{eq0201}) can be easily estimated for each of the families. In some cases, 
one of the families corresponds to the contribution of the precursor wave.

Note that the local contour deformation has been introduced in \cite{Shanin}. 
In \cite{Miklowitz} and \cite{Randles} the authors use a global contour transformation in order to get a comprehensive description of the transient phenomena.  

{\bf 6.} The deformation of the contour described above makes sense if the structure of the 
integral becomes more clear. This happens when the integrand decreases on all important branches as the result of deformation. Since the sum or the integrals over the 
branches remains the same, 
this means that some interference of the modes becomes eliminated.

Consider the exponential factor of (\ref{eq0201}): 
\[
\exp \{ i k_n (\omega) L - i \omega t \}.
\]   
Fix $L$ and $t$. Find under which condition this exponential factor decays when the 
(complex) point $\omega$ is shifted into the upper half-plane.  
Note that 
\begin{equation}
\frac{\ptl}{\ptl {\rm Im}[\omega]} \log (|\exp \{ i k_n (\omega) L - i \omega t \}|) = 
\left( t - \frac{\ptl {\rm Im}[k_n]}{\ptl {\rm Im}[\omega]} L \right) .
\label{eq0206}
\end{equation}
Due to Cauchy--Riemann conditions, 
\[
\frac{\ptl {\rm Im}[k_n]}{\ptl {\rm Im}[\omega]}
=
\frac{\ptl {\rm Re}[k_n]}{\ptl {\rm Re}[\omega]}
\equiv v_{{\rm gr},n}^{-1}(\omega),
\]
where $v_{{\rm gr},n}$ is the group velocity of given branch of the 
dispersion diagram at the given frequency.
Thus, the contour can be shifted into the upper half-plane if the value in the right is negative, i.~e.\ if 
\begin{equation}
v_{{\rm gr},n} < \frac{L}{t}.
\label{eq0207}
\end{equation}  
This relation explains why the method is efficient for finding the precursor 
waves. By definition, the precursor wave travels faster than any group velocity 
found on the 
{\em real\/} dispersion diagram for the given frequency range. 
Thus, for the values 
$L/t > {\rm max}(v_n)$ the contour can be deformed.   
As the result of the deformation, the branches change, and for complex $\omega$ some group velocities become bigger. These bigger velocities correspond to the precursors. 
Thus, the group velocity of the precursor wave is equal to 
\[
v_{\rm gr} = \left(\frac{\ptl {\rm Re}[k_n]}{\ptl {\rm Re}[\omega]} \right)^{-1}
\]
for some branch of the dispersion diagram for some complex $\omega$.

{\bf 7.}
The precursor waves, typically, are exponentially decaying ones. 
This is explained by the fact that they are close to leaky waves in multilayer media. 
It is easy to estimate the decay of a precursor wave. Again, consider the exponential factor. 
Note that   
\[
|\exp \{ i k_n (\omega) (L + \Delta L) - i \omega (t + \Delta t) \}| = 
\]
\[
\exp \{ - {\rm Im} [k (\omega)] \Delta L + {\rm Im} [\omega] \Delta t\}
\cdot
|\exp \{ i k_n (\omega) L  - i \omega t \}|.
\]
Since $\Delta t$ can be estimated as $L / v$, the precursor wave decays as $~e^{-\kappa L}$,
where 
\begin{equation}
\kappa = {\rm Im}[k_n(\omega)] - \frac{{\rm Im}[\omega]}{v},
\label{eq0208}
\end{equation}
where the complex value of $\omega$ and the branch $k_n$ correspond to the precursor.  

\subsection{Dispersion diagram and its analytical continuation}
\label{sec_overview}

\begin{figure}[ht]
\centerline{\epsfig{file=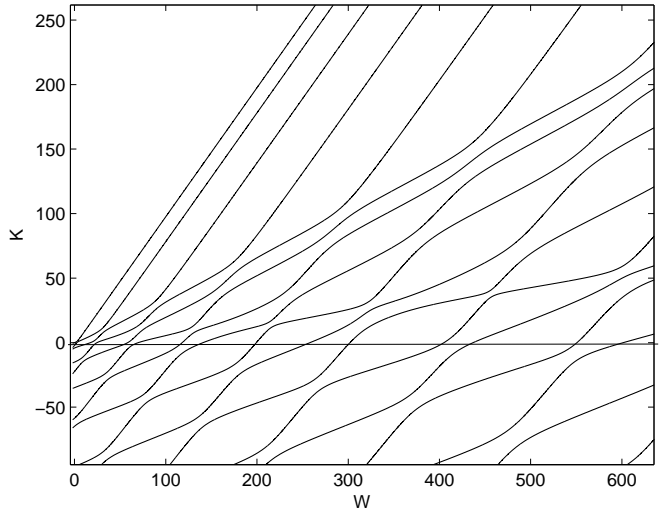, width = 10cm}}
\caption{An example of the dispersion diagram for a 3-layer waveguide. Real $W$} 
\label{fig04}
\end{figure}

Let us demonstrate some basic properties of dispersion diagrams of a 3-layer waveguide.
An example of the dispersion diagram is shown in Fig.~\ref{fig04}. The following parameters 
are used for numerical computations: 
\[
H_1 = 1, 
\qquad 
H_2 = 2, 
\qquad 
H_3 = 2.6, 
\]
\[
c_1 = 1, 
\qquad 
c_2 = 1.7 , 
\qquad 
c_3 = 3.2, 
\]
\[
\rho_1 = 15, 
\qquad 
\rho_2 =  1, 
\qquad 
\rho_3 =  1 .
\]
All values are assumed to be dimensionless. 

The dispersion diagram is represented in the coordinates $(W,K)$. These coordinates are convenient for the following reasons: a) they enable one to display propagating modes ($K > 0$)
and evanescent modes ($K < 0$) on a single real graph; b) for a single layer with height $h$, 
velocity $c$, and Newmann boundary condition the diagram has form of a family of parallel 
straight lines: 
\begin{equation}
K = \frac{W}{c^2} - \frac{\pi^2 n^2}{h^2}, 
\qquad
n = 0 , 1 , 2 ,\dots  
\label{eq0209}
\end{equation}

Fig.~\ref{fig04} represents the ``usual'' dispersion diagram, 
i.~e.\ $W$ is positive real. In the domain $W>0$, $K>0$ one can see clearly three zones 
occupied by the branches of the dispersion diagram: 
$W/c_2^2 <K < W/c_1^2$,
$W/c_3^2 <K < W/c_2^2$,
$0 <K < W/c_3^2$.
In the first zone one cans see, roughly, the modes corresponding to the layer of width $h_1$
with velocity $c_1$. In the second zone there is a terraced structure related to the waves with velocities $c_1$ and $c_2$. In the third zone the situation is less clear, but we can expect a complicated terraced structure caused by a mixture of modes of three types. 

\begin{figure}[ht]
\centerline{\epsfig{file=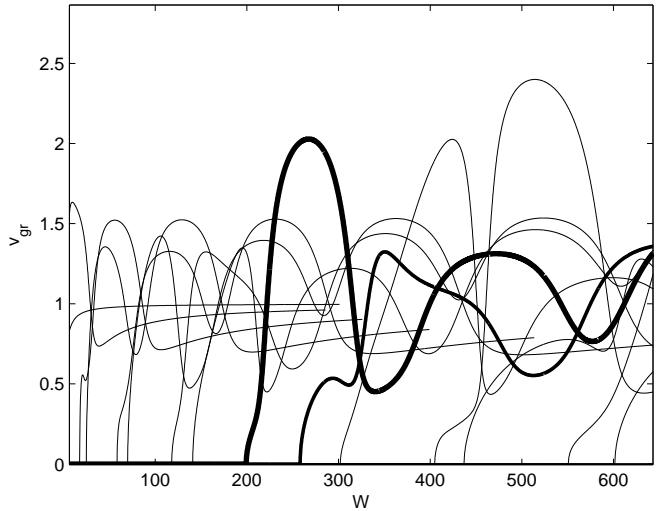, width = 10cm}}
\caption{Group velocities of the 3-layer dispersion diagram} 
\label{fig05}
\end{figure}

The group velocities of the modes are shown in Fig.~\ref{fig05}. The velocities are computed by the usual formula $v_{\rm gr} = (dk_n / d\omega)^{-1}$. Two random modes are made bold 
to demonstrate a complicated behavior of the curves.      
One can see that there are no modes having group velocity higher than $3$
within the frequency range displayed in the figure. 
Moreover, the group velocities are lower than~$2.5$ in the frequency range $0<W<300$. However, the velocity $c_3$ is equal to 3.2, thus, one can expect that there exists a wave traveling through the waveguide with this velocity. This wave can be observed in numerical modeling \cite{Shanin},  and it is a precursor. It decays exponentially in the waveguide. 
Since there are no physical losses, the decay of the precursor is caused by diffraction. 
The energy leaks into the guided waves traveling with smaller velocities. 

Let us study the analytical continuation of the dispersion diagram. 
Namely, for any given positive $\Omega$ consider the branches of $k(\omega' + i\Omega)$, 
where $\omega'$ is a real variable. With respect 
to the variable $W$, we study the curves
\[
W = (\sqrt{W'} + i \Omega)^2,
\] 
where $W'$ is real. 

\begin{figure}[ht]
\centerline{\epsfig{file=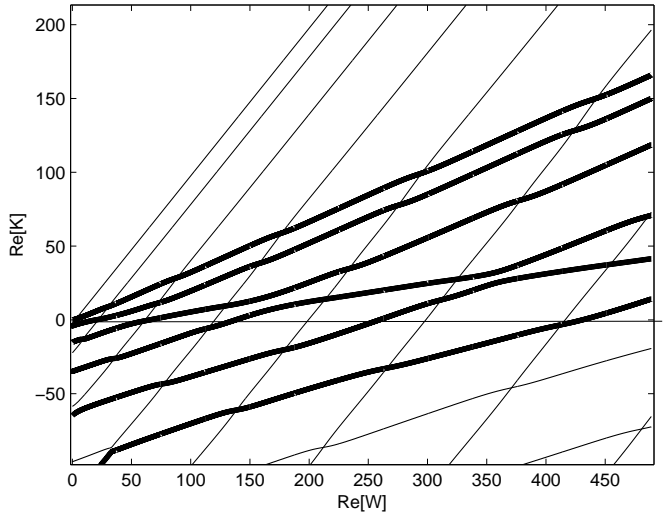, width = 10cm}}
\caption{Analytical continuation of the dispersion diagram for $\Omega=1$} 
\label{fig06}
\end{figure}

For the same physical configuration, and for $\Omega = 1$, the real part of the dispersion diagram is 
plotted in Fig.~\ref{fig06}.  One can see that the terrace-like structure partly disappears. 
Namely, the families of curves having slope $d K / d W $ close to $c_1^{-2}$ and $c_2^{-2}$
 now cross each other.
So do curves with slopes $c_1^{-2}$ and $c_3^{-2}$. However, the curves with slopes 
$c_2^{-2}$ and $c_3^{-2}$ still cannot cross and form so-called ``pseudo-crossings'' \cite{Shanin}.

Such a behavior can be interpreted from the physical and from the mathematical point of view. 
From the mathematical point of view, globally there are sheets corresponding to each layer separately (see below). These sheets are connected by branch points. The line $\omega = \omega' + i$ passes above the 
branch points connecting the sheets corresponding to layers 1 and 2. Also it passes above the 
branch points connecting the sheets corresponding to the layers 1 and 3. However, it passes below the branch points connecting the sheets corresponding to the layers 2 and 3. 
At the end of the paper we compute the positions of the branch points mentioned here 
and prove this conclusion.

Note that the crossings of the modes in Fig.~\ref{fig06} are just crossings of the real part of~$K$.
The imaginary parts of $K$ for crossing curves are usually different (otherwise such a crossing is a
branch point).  

The physical interpretation is as follows. The initial contour of integration (the real axis) and the 
diagram Fig.~\ref{fig04} correspond to a modal expansion of the field in a usual sense. Each mode is a complicated mixture of the fields in each of three layers. 

\begin{figure}[ht]
\centerline{\epsfig{file=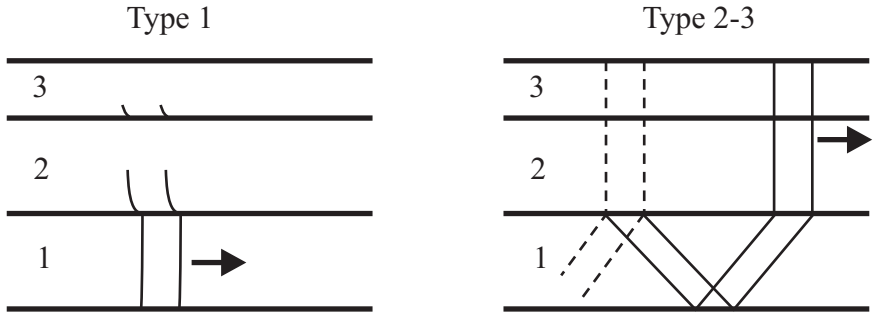}}
\caption{Structure of modes of Type~1 and Type~2-3} 
\label{fig07}
\end{figure}

The integration along the contour with $\Omega = 1$ corresponds to the expansion, in which two types of modes participate. They are modes of ``Type~1'' and of ``Type~2-3''. The modes of type 1 correspond to almost straight lines with $d {\rm Re}[K]/ d {\rm Re}[W] \approx c_1^{-2}$. 
The modes of Type~2-3 are all other curves (they are made bold in the figure). Schematically these modes are shown in Fig.~\ref{fig07}. Modes of Type~1 can be considered as surface waves traveling mainly in layer~1 and having exponentially decaying traces in the other layers. Modes of Type~2-3 have a complicated mixed 
structure across the layers 2 and 3, and they are leaky into the (slower) layer~1. 

\begin{figure}[ht]
\centerline{\epsfig{file=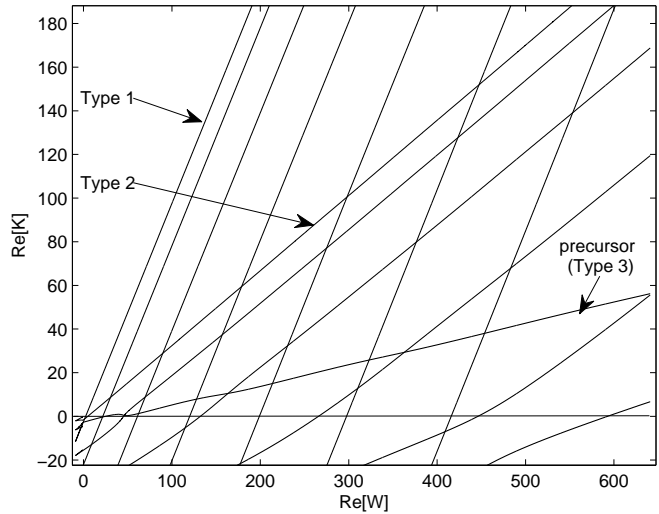, width = 10cm}}
\caption{Dispersion diagram for $\Omega = 3$ } 
\label{fig08}
\end{figure}

\begin{figure}[ht]
\centerline{\epsfig{file=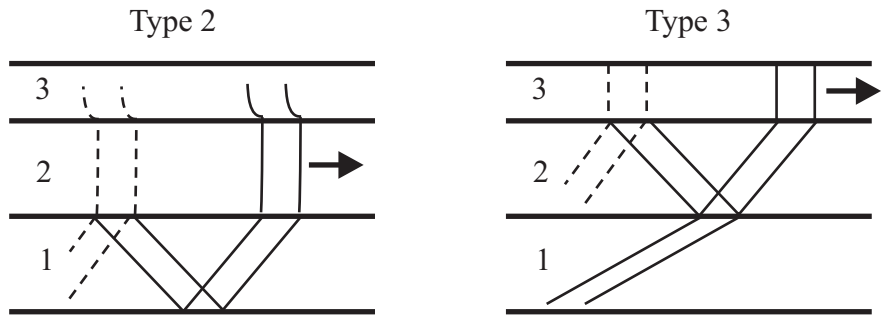}}
\caption{Structure of modes of Type~2 and Type~3} 
\label{fig09}
\end{figure}

Now consider the dispersion diagram on the contour with ${\rm Im}[\omega] = \Omega = 3$. The real part of the 
dispersion diagram is shown in Fig.~\ref{fig08}. One can see that all lines are now almost straight. 
Mathematically this means that the integration contour 
$\omega = \omega' + 3i$
passes above all branch points. 
Physically, the fields becomes expanded as series over modes of Type~1, Type~2, and Type~3. 
Type~1 has been explained above, and the sketches of Type~2 and Type~3 are shown in Fig.~\ref{fig09}. 
  
Modes of Type~2 are ``propagating'' in layer~2, ``exponentially decaying'' in the (faster) layer~3,
and ``leaky'' into the (slower) layer~1. Note that the the terms 
``propagating'', 
``exponentially decaying'', and 
``leaky'' have their exact meaning only for real~$\omega$. For complex $\omega$, say, for Type~2 mode, both the ``leaky'' part and the ``exponentially decaying'' part are 
decaying exponentially in the transversal direction.

There is a single branch in the figure corresponding to Type~3, it is the fastest wave in the system, and, thus, it is the first precursor in the waveguide. 

According to item~6 of the previous section (and according to a careful study of the group velocities from Fig.~\ref{fig04}, \ref{fig06}, \ref{fig08}), representation with $\Omega = 0$
can be used everywhere, representation with $\Omega = 1$ is useful for 
$L > 2.5t$, and the representation with $\Omega = 3$ is useful to describe the fastest precursor moving with the velocity $L/t \approx 3$.

To conclude this section, let us describe briefly the procedure of building of analytical continuations of the dispersion diagram. The procedure comprises two steps.    
On the first step a reference point $\omega= \omega_0$ is taken, and the values
$k_n(\omega_0)$ are found approximately. We take positive imaginary $\omega_0$, thus
all $k_n(\omega_0)$ are expected to be positive imaginary as well. They can be found by plotting a graph $D (\omega , k)$ for positive imaginary $k$ and finding the roots graphically (i.~e.\ as the points where the graph crosses the line $D=0$). Then, a contour in the plane $\omega$ is set. It should 
start at the point $\omega_0$ and then go along the line ${\rm Im}[\omega] = \Omega$. 
The contour is meshed, i.~e.\ a dense set of nodes $\omega_n$ is put on it.   
We move from node to node (starting at $\omega_0$) and find the roots of equation $D(\omega_n, k) =0$ by the Newton's method (by iterations). At each step we use the roots found on the previous steps as the initial approximations.  


\subsection{Demonstration of numerical capabilities of Miklowitz--Randles method}
\label{sec_numerical}

To demonstrate the capabilities of the Miklowitz--Randles method described above we consider a sample 
problem. The waveguide  is excited 
by the inhomogeniety (\ref{eq0104}). The point of excitation is located on the upper surface:
$y_0 = H_3$.    
The excitation  problem can be solved as follows. 
The field 
is represented as a modal expansion: 
\begin{equation}
\tilde u (x, y, t) = \sum_n \int \limits_0^{\infty} 
\hat f (\omega) P_n (\omega) U_n(\omega, y) 
\exp \{ \i k_n (\omega) |x| - \i \omega t \}   
d \omega
\label{eq0401}
\end{equation}
where $u$ is $u_1$, $u_2$, or $u_3$ (depending on $y$),  $U_n(\omega , y)$ is the mode profile:
\[
U_n (\omega , y) = 
\left\{ \begin{array}{ll} 
A_n(\omega) \cos(\alpha_1 (\omega, k_n) y) , & 
0< y < H_1, 
\\ 
B_n(\omega) \sin(\alpha_2 (\omega ,k_n) y) + C(\omega) \cos(\alpha_2 (\omega, k_n) y) , & 
H_1< y < H_2 , 
\\ 
F_n(\omega) \cos(\alpha_3 (\omega, k_n) y) , & 
H_2 < y < H_3, 
\end{array} \right. 
\]
The coefficients $A, B, C, F$ should obey the relation 
\begin{equation}
{\bf D}(\omega , k_n) \left( \begin{array}{c}
A_n(\omega) \\ B_n(\omega) \\ C_n(\omega) \\ F_n(\omega)
\end{array} \right) 
=
 \left( \begin{array}{c}
0 \\ 0 \\ 0 \\ 0
\end{array} \right).
\label{eq0402}
\end{equation}
Since ${\rm det}({\bf D}(\omega , k_n)) = 0$ and coefficients $A$ and $F$ cannot be equal to zero, 
the coefficients can be found as follows: 
\begin{equation}
F = 1, 
\qquad 
\left(  \begin{array}{c}
A_n(\omega) \\ B_n(\omega) \\ C_n(\omega)
\end{array} \right)
=
-
\left(  \begin{array}{ccc} 
d_{1,1} & d_{1,2} & d_{1,3} \\
d_{2,1} & d_{2,2} & d_{2,3} \\
d_{3,1} & d_{3,2} & d_{3,3} 
\end{array} \right)^{-1} 
\left(  \begin{array}{c}
d_{1,4} \\
d_{2,4} \\
d_{3,4}
\end{array} \right)
\label{eq0403}
\end{equation}
where $d_{m,n}$, $m,n = 1,\dots,4$ are elements of matrix ${\bf D} (\omega , k_n)$. 

The values $P_n(\omega)$ are amplitudes of the modes. To find them, note that 
the singularity in equation (\ref{eq0104}) is equivalent to a discontinuity of the $x$-derivative 
of $\tilde u$ at $x =0$. Namely, the following relation should be valid: 
\begin{equation}
i k_n (\omega) P_n(\omega) \sum_n U_n(\omega , y) = \frac{1}{2} \delta(y-y_0) .
\label{eq0404}
\end{equation}  
The values $P_n$ can be found by noting that the modes are orthogonal (see Appendix~A): 
\begin{equation}
<U_n , U_m> \equiv
\int \limits_0^{H_3} \rho(y) U_n (\omega,y) U_m(\omega,y) dy = 
0
\qquad
\mbox{ for }
n \ne m.
\label{eq0405}
\end{equation}
Thus, 
\begin{equation}
P_n(\omega) = \frac{1}{2 i k_n (\omega)} 
\frac{\rho(y_0) U_n (\omega, y_0)}{<U_n (\omega , y) , U_n (\omega , y)>}
\label{eq0406}
\end{equation}
Note that $A, B, C$, and thus $<U_n (\omega , y) , U_n (\omega , y)>$ can be found 
explicitly and they are trigonometric functions of~$\alpha_j(\omega, k_n)$. 

Since the observation point is taken on the upper surface ($y_0 = H_3$), and $F \equiv 1$, function 
$Y_n(\omega)$ from (\ref{eq0201}) is expressed as 
\[
Y(\omega, k_n (\omega)) = P_n(\omega) U_n(\omega , y_0)  = P_n(\omega).
\]
Note that $Y_n$ has the same Riemann surface as $k(\omega)$. 

For numerical demonstration of Miklowitz--Randles method we choose the 
spectrum of the excitation function 
\[
\hat f(\omega) = \exp \{ - (|\omega|-12)^2 / 16 \}. 
\]  
Such a spectrum is localized in the domain shown in Fig.~\ref{fig04}.
The distance from the source to the observation point $L$ is equal to
$L=10$. 

We compute $u_3(L,H_3,t)$ in three ways. First, we use formula (\ref{eq0201})
with integration along the real axis of $\omega$. The summation is performed 
over all modes shown in Fig.~\ref{fig04} (there are 19 of them). Of course, the series is truncated, 
so the result is not exact, but all propagating modes within the frequency band 
of the source
are taken into account. Therefore the error should be small. The result 
of the computation will be referred to as $u_{\rm a} (L , H_3 , t)$.
This signal will be used as a reference.  

Second, we use formula (\ref{eq0201}) performing integration along the contour passing 
mostly along the line ${\rm Im}[\omega] = 1$. The dispersion diagram on this line is shown in 
Fig.~\ref{fig06}. Only the modes of Type~2-3 are taken in the summation (i.~e.\ all modes of Type~1
are omitted). There are totally 6 modes that belong to Type~2-3. they are marked by bold lines 
in Fig.~\ref{fig06}. The result of the computations is referred to as $u_{\rm b} (L , H_3 , t)$.

Third, we use formula (\ref{eq0201}) performing integration along the contour passing 
mostly along the line ${\rm Im}[\omega] = 3$, i.~e.\ for the dispersion diagram shown in 
Fig.~\ref{fig08}. For the summation we take {\em just a single\/} mode corresponding to 
Type~3 (i.~e.\ no summation is held). The result is $u_{\rm c}(L , H_3 , t)$. 

\begin{figure}[ht]
\centerline{\epsfig{file=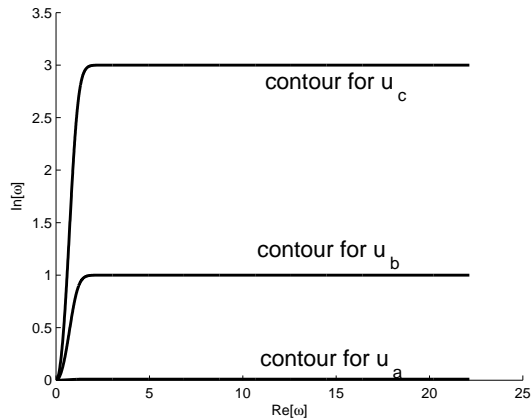, width = 8cm}}
\caption{Integration contours for $u_{\rm a}$, $u_{\rm b}$, $u_{\rm c}$} 
\label{fig12}
\end{figure}

The contours of integration are shown in Fig.~\ref{fig12}. 

\begin{figure}[ht]
\centerline{\epsfig{file=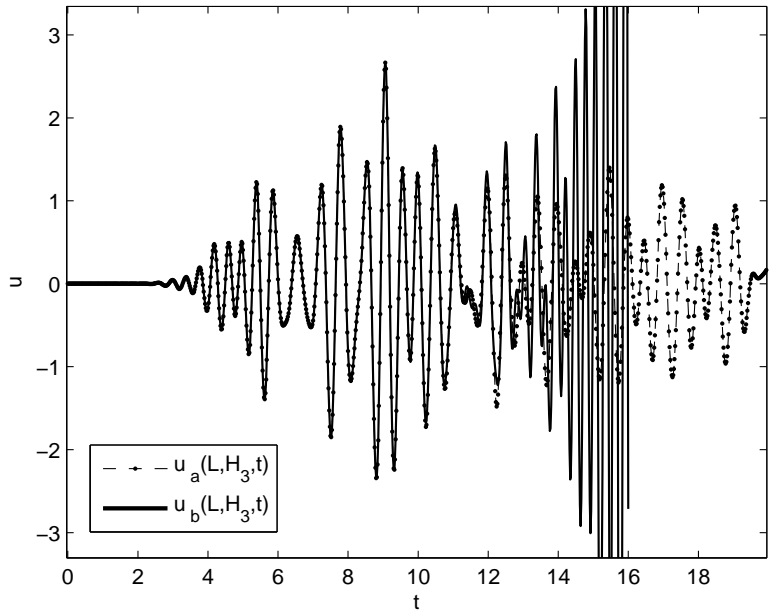, width = 12cm}}
\caption{Graphs of $u_{\rm a}$ and $u_{\rm b}$} 
\label{fig13}
\end{figure}

\begin{figure}[ht]
\centerline{\epsfig{file=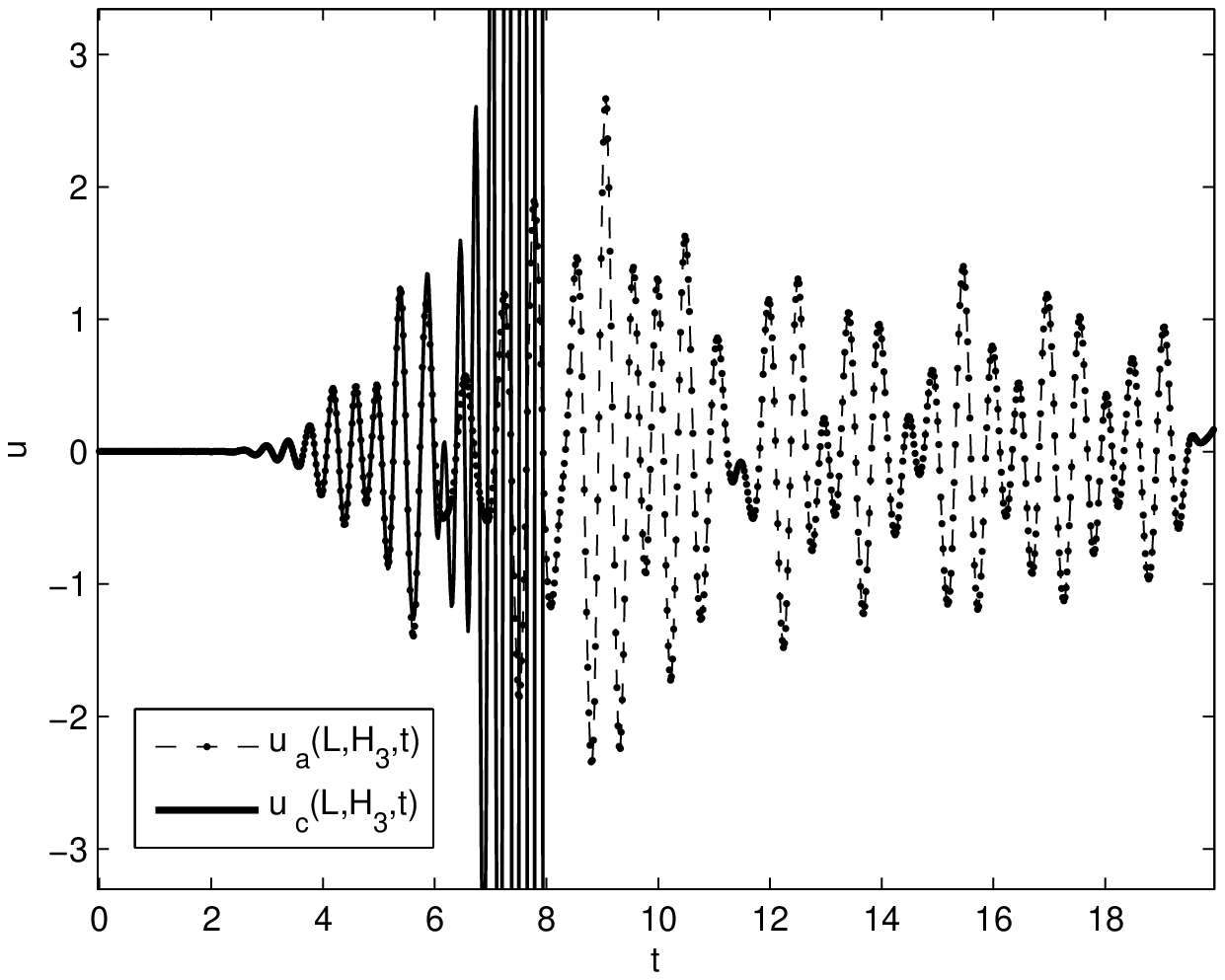, width = 12cm}}
\caption{Graphs of $u_{\rm a}$ and $u_{\rm c}$} 
\label{fig14}
\end{figure}

The results of computations are shown in Fig.~\ref{fig13} and Fig.~\ref{fig14}. 
In Fig.~\ref{fig13} the graphs of $u_{\rm a} (L , H_3 , t)$
and $u_{\rm b} (L , H_3 , t)$ are plot together. One can see that the graphs are close to each other 
for $t < 11$. This conclusion agrees with our prior summary of the Miklowitz--Randles
method. Namely, the omitted terms correspond to waves of Type~1 traveling with velocity 
$c_1 = 1$. Thus, we can expect a large error for $t > L / c_1$. 

In Fig.~\ref{fig14} the graphs of  $u_{\rm a} (L , H_3 , t)$
and $u_{\rm c} (L , H_3 , t)$ are displayed. One can see that the graphs are close to each other for 
$t < 6$. Again, since we omitted the waves of Type~1 and Type~2, we can expect a large error
for $t > L / c_2$.    

Figures \ref{fig13} and \ref{fig14} illustrate the main idea of the Miklowitz--Randles method. 
By deformation the integration contour in (\ref{eq0201}) and selection of a few branches 
corresponding to the fastest waves, one can achieve a good description of the fastest components 
of the wave field (i.~e.\ of the precursors).    
There are two types of the precursors: relatively slow waves of Type~2-3 and fast waves 
of Type~3. Both waves are faster than the guided waves in the far-field zone (the group velocities 
of the real dispersion diagram). 

\section{Finding the position of branch points of the dispersion diagram}
\label{sec_branchpoints}

\subsection{A family of waveguide problems with massive interfaces}

The structure of the Riemann surface of the dispersion diagram is important for understanding of physical processes in a waveguide.  
Here we develop an approach that can bring some order into this structure. Note that 
a similar problem for an elastic plate has been studied in \cite{Randles} using a different (explicit) approach. 

\begin{figure}[ht]
\centerline{\epsfig{file=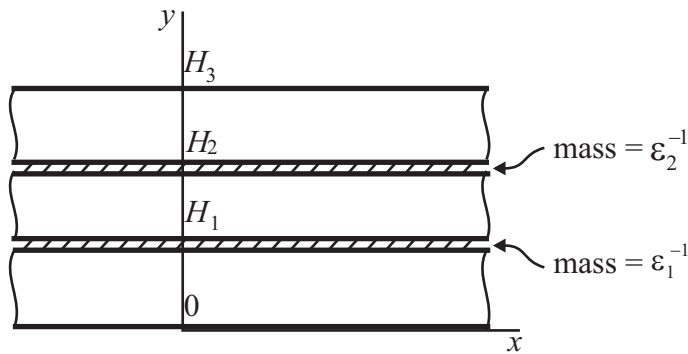}}
\caption{Additional interface masses} 
\label{fig01a}
\end{figure}

The idea is to consider a three-layer waveguide with modified interfaces between the layers.  
Namely, governing equations (\ref{eq0108}), 
boundary conditions (\ref{eq0109}), (\ref{eq0110}) remain the same, and boundary conditions 
(\ref{eq0111}) are rewritten as follows: 
\begin{equation}
\rho_2 u (H_1+0) - \rho_1 u (H_1-0) = \eps_1^{-1} u' (H_1 -0), 
\label{eq0301a}
\end{equation}
\begin{equation} 
\rho_3 u (H_2+0) - \rho_2 u (H_2-0) = \eps_2^{-1} u' (H_2+0), 
\label{eq0301b}
\end{equation}
where $\eps_1$ and $\eps_2$ are the {\em linking parameters}. 
One can understand (\ref{eq0301a}), (\ref{eq0301b})
as the introduction of some mass of the interfaces between the layers. 
Namely, recall that $u_j$ is the acoustical potential. 
(\ref{eq0301a}) can be rewritten as 
\[
p_2 - p_1 = -\eps_1^{-1} a,
\]
where $p_2$ is the acoustical pressure at $y = H_1 + 0$, $p_1$ is the acoustical 
pressure at $y = H_1 - 0$, $a$ is the $y$-component of the acceleration of the 
points of the boundary between the layers 1 and~2. Thus, $\eps_1^{-1}$ can 
be interpreted as the mass of the interface line per unit length. 
The positions of the interfaces are shown in Fig.~\ref{fig01a}. 

Now the values $k_n(\omega)$ depend implicitly on parameters $\eps_1$ and $\eps_2$:
\[
k_n = k_n(\omega , \eps_1 , \eps_2).
\]

\subsection{A method to trace the branch points}

Thus, we introduced a family of waveguide problems depending on parameters $\eps_1, \eps_2$. 
In this family we study two special cases: 
(a) the case $\eps_1 , \eps_2 \to 0$ corresponding to 
very heavy interfaces, i.~e.\ to
three disjoint layers with Neumann 
boundary conditions; 
(b) the case $\eps_1 , \eps_2 \to \infty$ corresponding
to lightweight interfaces, i.~e.\ 
to the problem formulated above (conditions (\ref{eq0301a}, (\ref{eq0301b})) tend to (\ref{eq0111})). 
Case (a) is relatively simple and can be studied by perturbation methods. Case (b) 
is our target. In the plane of parameters $(\eps_1 , \eps_2)$ case (a) corresponds to a point asymptotically close to the origin, while case (b) corresponds to a point close to the infinity. 
The key idea of this section is to connect the origin and the infinity by some path 
and study a continuous transition between case (a) and case (b). This 
can be physically interpreted as 
a gradual loss of mass of the interfaces. 

We prefer not to think about $\eps_{1,2}$ as about physical parameters. 
Instead, we allow $\eps_{1,2}$ to take arbitrary complex values. 
We assume that functions $k_n$ depend on $\eps_{1,2}$ (as on complex parameters)
{\em analytically\/} almost everywhere. This is true since the dispersion equation depends 
on $\eps_1$ and $\eps_2$ analytically (see below). This analyticity has an important corollary. Namely, the topology of the Riemann surface for case (b) 
is the same as for case (a) (to be more precise, there exists a homotopy between the surfaces (a) and (b)). Particularly, each branch point of the Riemann surface of the dispersion diagram for 
case (b) corresponds to a branch point for case~(a)
(and vice versa). The branch points in case (a) are ordered in a 
very clear way. 
Practically, the branch points 
are indexed by four indexes.
By continuation, we assign the same indexes to the branch points for case~(b). 

Note that the possibility to make $\eps_1$, $\eps_2$ complex is important only from the 
theoretical point of view. It is necessary to prove the analiticity of the solution 
with respect to these parameters. In the practical computations $\eps_1$, $\eps_2$
will be always real. 

We prefer to study the problem in the coordinates $(W,K)$ instead of $(\omega, k)$,
particularly, we introduce 
\[
K_n (W) \equiv (k_n (\sqrt{W}))^2.
\] 

The dispersion equation for the system with interfaces of variable mass
is as follows: 
\begin{equation}
\tilde D(W, K , \eps_1 , \eps_2) = {\rm det} ({\bf D}_0 + \eps_1 {\bf D}_1 + \eps_2 {\bf D}_2) = 0 , 
\label{eq0302}
\end{equation} 
\begin{equation}
{\bf D}_0 =
\left( \begin{array}{cccc}
 \alpha_1 \sin (\alpha_1 H_1) & 0 & 0 & 0  \\
\alpha_1 \sin (\alpha_1 H_1) & \alpha_2 \cos(\alpha_2 H_1) & -\alpha_2 \sin (\alpha_2 H_1) & 0  \\
0 & 0 & 0 &  -\alpha_3 \sin (\alpha_3 h_3) \\
0 & \alpha_2 \cos (\alpha_2 H_2) & 
-\alpha_2 \sin(\alpha_2 H_2) &  -\alpha_3 \sin (\alpha_3 h_3) 
\end{array} \right) ,
\label{eq0303}
\end{equation} 
\begin{equation}
{\bf D}_1 =
\left( \begin{array}{cccc}
-\rho_1 \cos (\alpha_1 H_1)  & \rho_2 \sin(\alpha_2 H_1) & \rho_2 \cos (\alpha_2 H_1) & 0  \\
0 & 0 & 0 & 0  \\
0 & 0 & 0 & 0 \\
0 & 0 & 0 & 0 
\end{array} \right) ,
\label{eq0304}
\end{equation} 
\begin{equation}
{\bf D}_2 =
\left( \begin{array}{cccc}
0 & 0 & 0 & 0  \\
0 & 0 & 0 & 0  \\
0 & -\rho_2 \sin (\alpha_2 H_2) & -\rho_2 \cos(\alpha_2 H_2) & \rho_3 \cos (\alpha_3 h_3) \\
0 & 0 & 0 & 0 
\end{array} \right) ,
\label{eq0305}
\end{equation} 

Consider case (a). Let $\eps_{1,2}$ be small positive parameters. In any finite domain of variables 
$(W,K)$ the dispersion diagram is represented almost everywhere by three families of branches: 
\begin{equation}
K_{\nu,n} (W) \approx \frac{W}{c_{\nu}^2} - \frac{\pi^2 n^2}{h_{\nu}^2} ,
\label{eq0306}
\end{equation}
where
$\nu = 1, 2, 3$ denotes the family, $n = 0 , 1, 2, \dots$ is the index of the branch in the 
family.
Obviously, these families of branches correspond to modes in separate layers with Neumann boundary conditions. 

There can be no branch points of the function $K(W)$ except in the vicinities of
the crossing points of the branches. These crossing points are denoted by 
$(W_{\mu, \nu , m,n} , K_{\mu, \nu, m,n})$.
(four indexes and no argument).
This notation corresponds to a crossing of the branches $K_{\mu , m}(W)$ and $K_{\nu ,n}(W)$, 
$\mu \ne \nu$. The positions of the crossings can be easily calculated from 
(\ref{eq0306}):
\begin{equation}
W_{\mu,\nu , m,n} = \pi^2 
\left( \frac{m^2}{h_{\mu}^2} - \frac{n^2}{h_{\nu}^2} \right) 
\left( \frac{1}{c_{\mu}^2} - \frac{1}{c_\nu^2} \right)^{-1},
\label{eq0307}
\end{equation}
\begin{equation}
K_{\mu,\nu , m,n} = \pi^2 
\left( \frac{m^2 c_{\mu}^2}{h_{\mu}^2} - \frac{n^2 c_{\nu}^2}{h_{\nu}^2} \right) 
\left( c_{\nu}^2 - c_\mu^2 \right)^{-1},
\label{eq0308}
\end{equation}

As it is shown in Appendix~B, generally, in the vicinity of a crossing there exists 
a pair of branch points of order~2 connecting corresponding sheets. I.~e.,\ near the 
crossing point $W_{\mu,\nu , m,n}$ there are two branch points 
connecting branches $K_{\mu , m}(W)$ and $K_{\nu , n}(W)$. 
Since $\tilde D(W,K,\eps_1,\eps_2)$ is a real function for real arguments, the branch points in the pair are complex conjugate to each other, i.~e.\ they can be written as 
\[
(W , K) = (\Theta_{\mu,\nu , m,n} (\eps_1 , \eps_2), 
\Xi_{\mu,\nu , m,n}(\eps_1 , \eps_2))
\]
\[
(W , K) = (\bar \Theta_{\mu,\nu,m,n} (\eps_1 , \eps_2),
\bar \Xi_{\mu,\nu,m,n} (\eps_1 , \eps_2)),
\]
where the bar denotes complex conjugation.  
These values are solutions of the system of equations following from (\ref{eq0204}):
\begin{equation}
\tilde D (W,K, \eps_1 , \eps_2)  = 0 , 
\label{eq0310}
\end{equation}
\begin{equation}
\tilde D_K (W,K, \eps_1 , \eps_2)  = 0 , 
\label{eq0311}
\end{equation}
where $\tilde D_K$ is the partial derivative of $\tilde D(W,K,\eps_1 , \eps_2)$ 
with respect to the second argument:
\[
\tilde D_K = \ptl_K \tilde D. 
\]

The crossing points $(W_{\mu,\nu,m,n} , K_{\mu,\nu,m,n})$ are limits of the positions
of the branch points:
\[
\lim_{\eps_1 , \eps_2 \to 0}
\Theta_{\mu,\nu,m,n} (\eps_1 , \eps_2)  = W_{\mu,\nu , m,n},
\qquad
\lim_{\eps_1 , \eps_2 \to 0}
\Xi_{\mu,\nu , m,n} (\eps_1 , \eps_2)  = K_{\mu,\nu, m,n}.
\]

\begin{figure}[ht]
\centerline{\epsfig{file=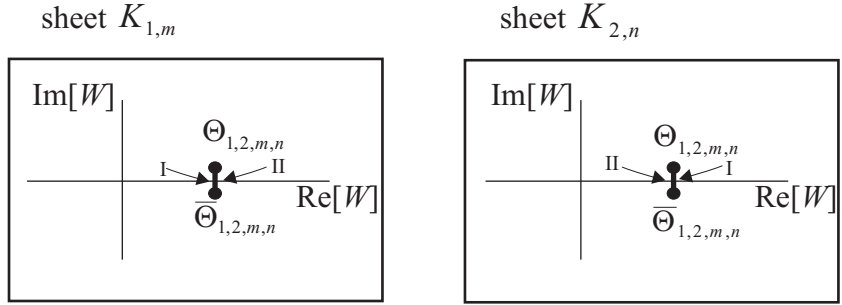}}
\caption{Local scheme of connection of sheets of Riemann surface} 
\label{fig10}
\end{figure}

Locally, an example of the scheme of two sheets of the Riemann surface of $K(W)$ is shown in 
Fig.~\ref{fig10}. 
A small cut is made between the points $\Theta_{\mu,\nu , m,n} (\eps_1 , \eps_2)$
and $\bar \Theta_{\mu, \nu , m,n} (\eps_1 , \eps_2)$ on the sheets 
$K_{\mu, m}$  and $K_{\nu , n}$. The shores of the cuts are denoted by Roman numbers I and~II.
The shores denoted by the same numbers are attached to each other.    
Thus, a connection between sheets exists, and it is local.                   
        
\begin{figure}[ht]
\centerline{\epsfig{file=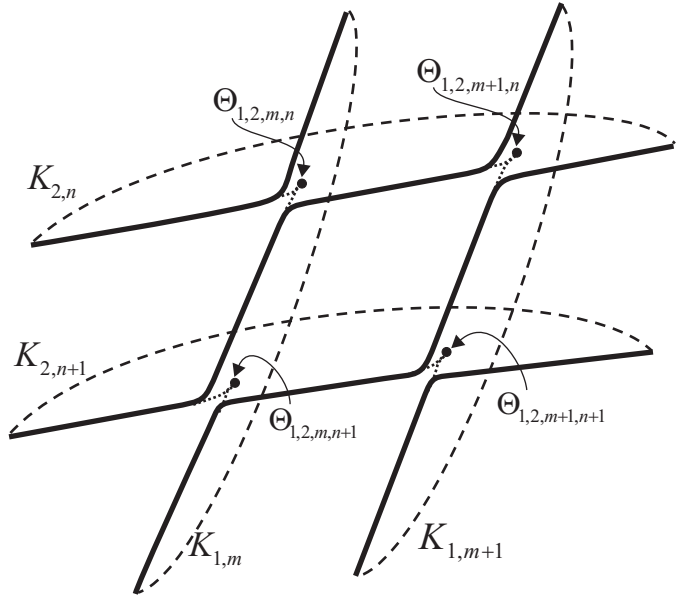}}
\caption{Formation of a terrace-like structure} 
\label{fig11}
\end{figure}

The nature of the terraced structure (at least for small $\eps_1$, $\eps_2$) becomes
clear. Since the real axis of $W$ passes between the points of each pair 
$(\Theta_{\mu,\nu,m,n}, \bar \Theta_{\mu,\nu,m,n})$,
each such pass leads to a transition from one sheet to another.  
In Fig.~\ref{fig11} we demonstrate several sheets of the Riemann surface cut along the real axis.
 
This understanding enables one to find numerically the the position of the 
branch points for case (b). Fix some family indexes $\mu,\nu$
(which can be equal to 1,2,3 and should be $\mu \ne \nu$) and some sheet indexes $m,n$
taking values $0,1,2\dots$. Consider the position of the branch point 
$\Theta = \Theta_{\mu,\nu, m, n}$ as a function of variables $\eps_1 , \eps_2$.
The target is to find the values 
\[
\lim_{\eps_1, \eps_2 \to \infty} \Theta (\eps_1 , \eps_2).
\]

\begin{figure}[ht]
\centerline{\epsfig{file=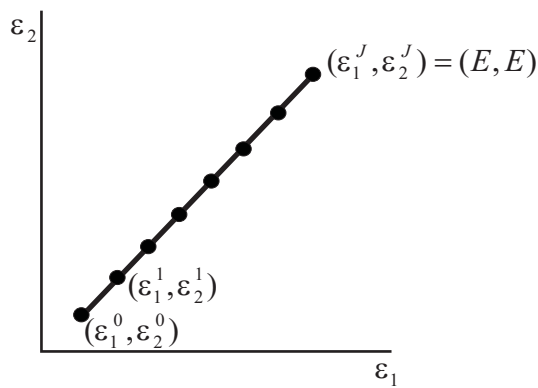}}
\caption{Contour in the $(\eps_1 , \eps_2)$ plane and its discretization} 
\label{fig16}
\end{figure}

The step-by-step procedure of finding 
these values is as follows. 

\noindent
--- Take some big real value $E$ and some small starting values $\eps_1^0 , \eps_2^0$
(0 is the upper index, not a power).
In the plane $(\eps_1 , \eps_2)$ draw a contour 
connecting the points $(\eps_1^0 , \eps_2^0)$ and $(E,E)$.
Put the nodes $(\eps_1^j , \eps_2^j)$ on this contour densely enough (see 
Fig.~\ref{fig16}). Let be 
$j= 1,\dots , J$, $\eps_{1,2}^J = E$. 

\noindent 
--- Using the formulas from Appendix~B, find the starting value 
$\Theta(\eps_1^0 , \eps_2^0)$. 
Find also the value of $K = K (\eps_1^0 , \eps_2^0)$ corresponding to this $\Theta$. 
Note that the process cannot be started with the values 
$\Theta(0 , 0) = W_{\mu, \nu , m,n}$, $\Xi(0 , 0)= K_{\mu,\nu , m,n}$, since the branch points collocate for these parameters, 
and the matrix ${\bf Q}$ (see below) is singular.  

\noindent
--- 
For $j = 1, \dots, J$ take $\eps_1 = \eps_1^j$, $\eps_2 = \eps_2^j$ and solve the system 
(\ref{eq0310}), (\ref{eq0311}) numerically with respect to $(W, K)$.
Use the Newton's method for this. 
Denote the solution of the system (\ref{eq0310}), (\ref{eq0311})
by $W = \Theta(\eps_1^j , \eps_1^j ) $, 
$K = \Xi(\eps_1^j , \eps_1^j )$.
As the starting approximation for the Newton's 
method, use the values $\Theta (\eps_1^{j-1} , \eps_1^{j-1})$ and 
$\Xi (\eps_1^{j-1} , \eps_1^{j-1})$ found on the previous step. 

--- After the last step, solve numerically the system 
\begin{equation}
\tilde D (W,K)  = 0 , 
\label{eq0312}
\end{equation}
\begin{equation}
\tilde D_K (W,K)  = 0 , 
\label{eq0313}
\end{equation}
corresponding to $\eps_1 = \eps_2 = \infty$ using 
$\Theta (\eps_1^{J} , \eps_1^{J})$ and 
$\Xi (\eps_1^{J} , \eps_1^{J})$ as the starting approximation. 
As the result, get $W = \Theta_{m,n,\mu, \nu}$, 
$K = \Xi_{m,n,\mu, \nu}$, which is the position of the 
branch point for the initial problem. 

The Newton's method used for solving the system (\ref{eq0310}), (\ref{eq0311})
for any fixed $\eps_1 , \eps_2$
or the system (\ref{eq0312}), (\ref{eq0313}) is as follows. The solution $(W, K)$ is found
by iterations. 
Some starting approximation is taken.
Then on each step the correction is made:
\begin{equation}
\left( \begin{array}{c}
W\\
K
\end{array} \right) 
\to  
\left( \begin{array}{c}
W \\
K
\end{array} \right)
+
\left( \begin{array}{c}
\Delta W\\
\Delta K
\end{array} \right) 
\label{eq0314}
\end{equation} 
\begin{equation}
\left( \begin{array}{c}
\Delta W\\
\Delta K
\end{array} \right)
= 
-
{\bf Q}^{-1}
\left( \begin{array}{c}
\tilde D(W, K) \\
\tilde D_K(W, K)
\end{array} \right),
\label{eq0314}
\end{equation} 
\[
{\bf Q} = {\bf Q} (W,K,\eps_1 , \eps_2)= 
\left( \begin{array}{cc}
\ptl_W \tilde D   & \ptl_K \tilde D  \\
\ptl_W \tilde D_K & \ptl_K \tilde D_K  
\end{array}\right). 
\] 
Obviously, 
(\ref{eq0314}) is a linear correction of the residue $(\tilde D(W,K), \tilde D_K (W, K))$.
In (\ref{eq0314}) the parameters $\eps_1 , \eps_2$ are assumed to be fixed.

We should comment the choice of the starting point $(\eps_1^0 , \eps_2^0)$. 
An obvious choice is to take $\eps_1^0 = \eps_2^0 = \eps$ for some small~$\eps$. 
This choice is good for the branch points $\Theta_{1,3,m,n}$ describing interaction between layers 1 and~3. For the branch points $\Theta_{1,2,m,n}$, however,
one should better take the initial values $\eps_1^{0} = \eps$, $\eps_2^{0} = 0$ (only layers 1 and 2 interact). 
Similarly, for the points $\Theta_{2,3,m,n}$ one should take 
$\eps_1^{0} = 0 $, $\eps_2^{0} = \eps$.

\subsection{Numerical results for tracing the branch points}

\begin{figure}[ht]
\centerline{\epsfig{file=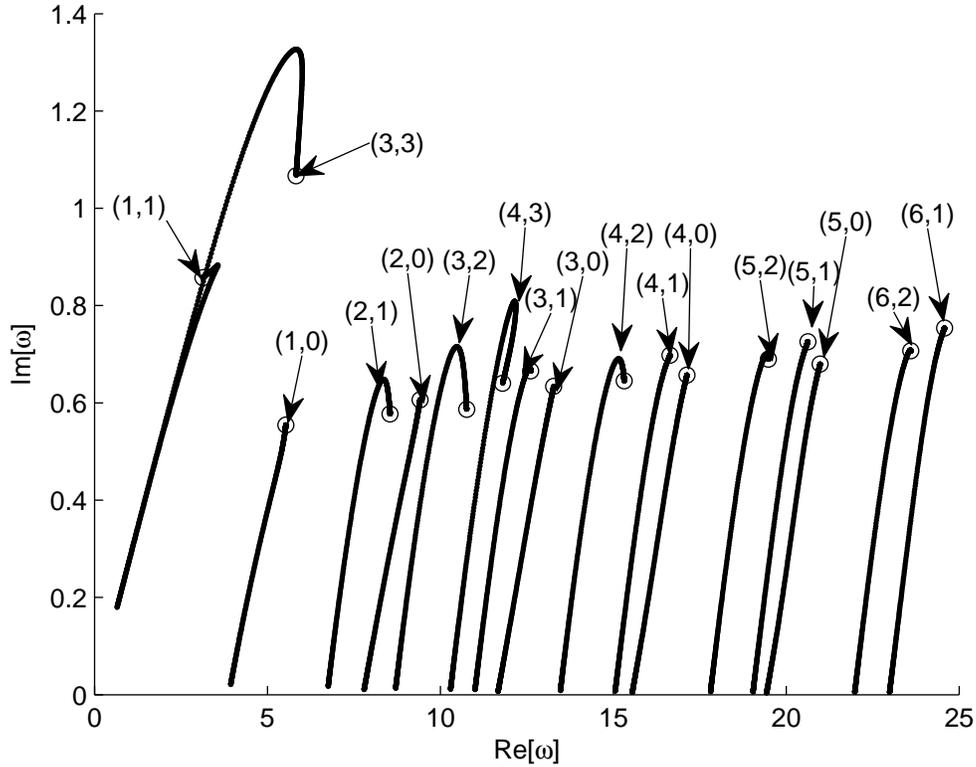}}
\caption{Tracing of branch points $\omega_* = \sqrt{\Theta_{1,2,m,n}}$. The pair $(m,n)$ is 
indicated near each trajectory} 
\label{fig_points_12}
\end{figure}

\begin{figure}[ht]
\centerline{\epsfig{file=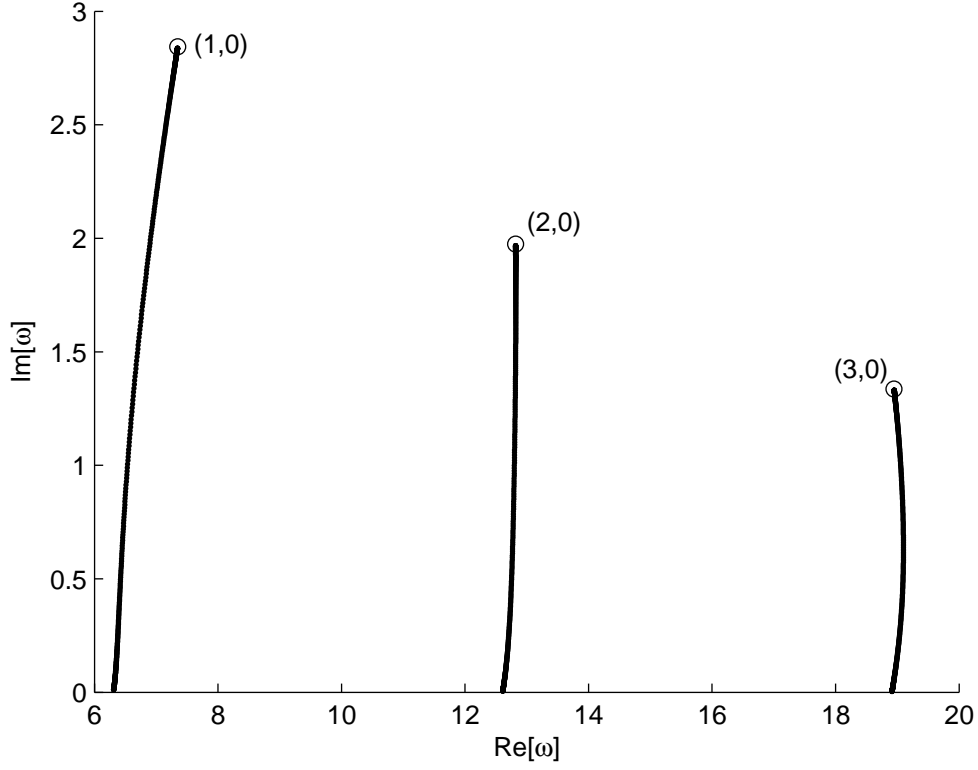}}
\caption{Tracing of branch points $\omega_* = \sqrt{\Theta_{2,3,m,n}}$. The pair $(m,n)$ is 
indicated near each trajectory} 
\label{fig_points_23}
\end{figure}

\begin{figure}[ht]
\centerline{\epsfig{file=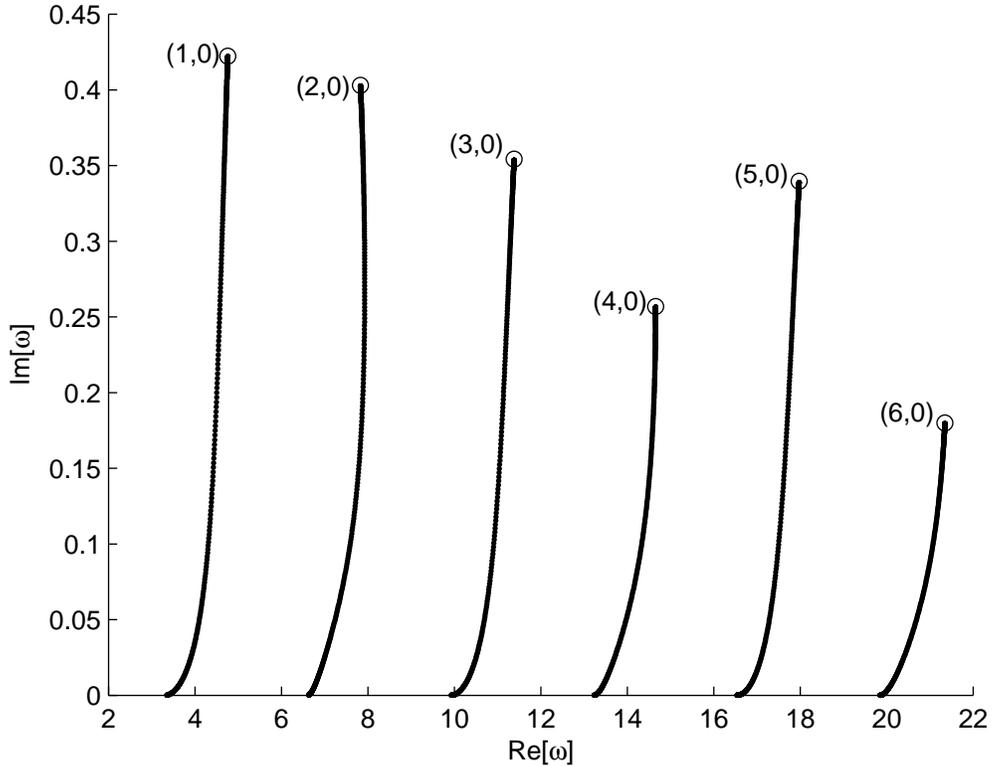}}
\caption{Tracing of branch points $\omega_* = \sqrt{\Theta_{1,3,m,n}}$. The pair $(m,n)$ is 
indicated near each trajectory} 
\label{fig_points_13}
\end{figure}

The procedure  of tracing the branch points has been implemented. 
The results are shown in Fig.~\ref{fig_points_12}, Fig.~\ref{fig_points_23}, and 
Fig.~\ref{fig_points_13}.  The positions of the 
branch points $\omega_* = \sqrt{\Theta_{\mu,\nu,m,n}}$ are shown everywhere.

In Fig.~\ref{fig_points_12}
the trajectories of branch points 
for $(\mu,\nu) = (1,2)$
traced along the contour 
\[
(\eps_1, \eps_2) = (\eps , \eps - 0.01),
\qquad  
\eps = 0.01 \, \dots \, 1000
\]
are shown.
The pair of indexes $(m,n)$ is plot near each trajectory. 
The last point is obtained by solving the system (\ref{eq0312}), (\ref{eq0313}),
and the result is plot as a small circle. 
We computed the positions of the branch points relevant for the  dispersion diagrams 
shown in Fig.~\ref{fig04}, Fig.~\ref{fig06}, and Fig.~\ref{fig08}. 

Each trajectory starts from vicinity of the real axis of $\omega$. This corresponds 
approximately to the 
points $\sqrt{W_{1,2,m,n}}$. In more details, we used 
the formulas (\ref{eqB19}), (\ref{eqB22}) to compute the starting positions of the 
branch points for $\eps_1 = 0.01$, $\eps_2 = 0$.     

The trajectories of the branch points for $(\mu,\nu) = (2,3)$ are shown 
in Fig.~\ref{fig_points_23}. The notations are the same. the starting values of the tracing are
$\eps_1 =0$, $\eps_2 = 0.01$.

Finally, the  trajectories of the branch points for $(\mu,\nu) = (1,3)$ are shown 
in Fig.~\ref{fig_points_13}. These points describe an interaction between the layers 
that are not neighbors, so complicated formulas (\ref{eqB38}), (\ref{eqB39}) are used for the first step. The values $\eps_1^0 = \eps_2^0 = 0.01$ are taken for computations.  

By comparing Figures~\ref{fig_points_12},~\ref{fig_points_23},~\ref{fig_points_13}
with Fig.~\ref{fig12} one can see that the points $\sqrt{\Theta_{2,3,1,0}}$,
$\sqrt{\Theta_{2,3,2,0}}$, $\sqrt{\Theta_{2,3,3,0}}$, $\sqrt{\Theta_{1,2,3,3}}$
are located between contours for $u_{\rm b}$ and $u_{\rm c}$, while all 
other branch points are located between the contours for $u_{\rm a}$ and $u_{\rm c}$. 
This conclusion is supported by analysis of the dispersion diagram Fig.~\ref{fig06}.
One can clearly see pseudo-crossings corresponding to the points
$\Theta_{2,3,2,0}$, $\Theta_{2,3,3,0}$.

\section{Conclusion}

Miklowitz--Randles' method is applied to the problem of wave propagation in a three-layer waveguide. It is demonstrated numerically that fast components of the signal can be computed 
as a  sum over a small subset of branches of the analytical continuation of the dispersion diagram.
Thus, the analytical continuation of the dispersion diagram is an important tool of analyzing the 
transient phenomena in layered waveguides. 

A physically motivated method is proposed to analyze the structure of the Riemann surface of the 
dispersion diagram. The waveguide is split into several layers linked by interfaces of variable 
linking parameters~$\eps_{1,2}$. The dispersion diagram depends on $\eps_{1,2}$ analytically, 
thus the structure of the Riemann surface is not changed as $\eps_{1,2}$ are varied. The case of 
small $\eps_{1,2}$ can be considered by using the perturbation method.

The positions of the branch points of the dispersion diagram of the layered waveguide can be 
obtained from the positions of the branch points for the auxiliary problem with small $\eps_{1,2}$
by applying a numerical tracing procedure. The validity of the method is demonstrated.   

The work is supported by the RSF grant 14-22-00042.

\section*{Appendix~A. Bilinear relations}

\subsection*{A1. Orthogonality relation for the modes of the initial problem}

Equation (\ref{eq0108}) can be written as follows: 
\begin{equation}
u'' + \frac{W}{c^2}  u - K u =0 
\label{eqC01}
\end{equation}
(the prime denotes the $y$-derivative).
Let there be two different solutions, $u_1(y)$ and $u_2(y)$, having common value 
of $W$ and different values of $K$, namely $K_1$ and $K_2$. Thus, 
\[
u_1'' + \frac{W}{c^2}  u_1 - K_1 u_1 =0,
\] 
\[
u_2'' + \frac{W}{c^2}  u_2 - K_2 u_2 =0,
\]
Multiply the first equation by $\rho(y) u_2$ and integrate (by parts) over the  
segments $[0, H_1]$, $[H_1, H_2]$, $[H_2 , H_3]$. Multiply the second equation
by $\rho(y) u_1$ and integrate over the same sum of the segments. 
Subtract the second result from the first result:
\[
\left . \rho_1 (u_1' u_2 - u_1 u_2') \right|_{0}^{H_1-0}
+
\left . \rho_2 (u_1' u_2 - u_1 u_2') \right|_{H_1+0}^{H_2-0}
+
\left . \rho_3 (u_1' u_2 - u_1 u_2') \right|_{H_2+0}^{H_3}
=
\]
\begin{equation}
(K_1 - K_2) \int \limits_0^{H_3} \rho(y) u_1(y) u_2(y) dy .
\label{eqC02}
\end{equation}
Due to boundary conditions (\ref{eq0109}), (\ref{eq0110}), (\ref{eq0111}) the left-hand side 
of this equation is equal to zero, and thus the right-hand side should be zero. 
If $K_1 \ne K_2$  then 
\begin{equation}
<u_1 , u_2> \equiv  \int \limits_0^{H_3} \rho(y) u_1(y) u_2(y) dy =0 .
\label{eqC03}
\end{equation}


\subsection*{A2. Bilinear relations for the problem with heavy interfaces}

Consider the problem comprised by equations (\ref{eqC01}), 
boundary conditions (\ref{eq0109}), (\ref{eq0110}), (\ref{eq0301a}), 
(\ref{eq0301b}). 
Let $u_1$ be a solution of this problem with parameters $W = W_1$, 
$K = K_1$, $\eps_1 = \eps_1^\dag$, $\eps_2 = \eps_2^\dag$  
for some arbitrary values. 
Let $u_2$ be a solution of this problem with parameters $W = W_2$, 
$K = K_2$, $\eps_1 = \eps_1^\ddag$, $\eps_2 = \eps_2^\ddag$.

Perform the same procedure as in the previous subsection (multiply the equation for 
$u_1$ by $\rho u_2$, multiply the equation for $u_2$ by $\rho u_1$, integrate over the layers,
and subtract the results). Using the boundary conditions 
(\ref{eq0109}), (\ref{eq0110}) get 
\begin{equation}
S +(W_1 - W_2) \int \limits_0^{H_3} \frac{\rho(y)}{c^2(y)} u_1(y) u_2(y) dy - 
(K_1 - K_2) \int \limits_0^{H_3} \rho(y) u_1(y) u_2(y) dy = 0 
\label{eqC04}
\end{equation} 
where 
\[
S = 
u_1'(H_1) (\rho_1 u_2(H_1-0) - \rho_2 u_2(H_1+0)) 
+
u_2'(H_1) (\rho_2 u_1(H_1+0) - \rho_1 u_1(H_1-0))
+
\]
\begin{equation}
u_1'(H_2) (\rho_2 u_2(H_2-0) - \rho_3 u_2(H_2+0)) 
+
u_2'(H_2) (\rho_3 u_1(H_2+0) - \rho_2 u_1(H_2-0))
\label{eqC05}
\end{equation} 
Consider some particular cases.
If $\eps_1^\dag$, $\eps_2^\dag$, $\eps_1^\ddag$, $\eps_2^\ddag$ are not equal to zero, then 
due to (\ref{eq0301a}), (\ref{eq0301b})   
\begin{equation}
S = 
[(\eps_1^\dag)^{-1} - (\eps_1^\ddag)^{-1}] u_1'(H_1) u_2'(H_1) 
+
[(\eps_2^\dag)^{-1} - (\eps_2^\ddag)^{-1}] u_1'(H_2) u_2'(H_2) 
\label{eqC06}
\end{equation}
Particularly, if $u_1$ and $u_2$ correspond to the same waveguide:
\[
\eps_1^\dag = \eps_1^\ddag , 
\qquad 
\eps_2^\dag = \eps_2^\ddag ,
\]
then $S = 0$. If $\eps_1^\dag = 0$ and $\eps_2^\dag =0 $ then 
$u_1' (H_1) = u_1'(H_2) =0$, and 
\[
S = 
\eps_1^2 [\rho_2 u_2(H_1+0) - \rho_1 u_2(H_1-0)]
[\rho_2 u_1(H_1+0) - \rho_1 u_1(H_1-0)]
+
\]
\begin{equation}
\eps_2^2
[\rho_3 u_2(H_2+0) - \rho_2 u_2(H_2-0)]
[\rho_3 u_1(H_2+0) - \rho_2 u_1(H_2-0)]
\label{eqC07}
\end{equation}

\subsection*{A3. Orthogonality for solutions at the branch points of the dispersion curve}

Fix some values $\eps_1$, $\eps_2$ and consider the problem composed of equation
(\ref{eqC01}) and boundary conditions (\ref{eq0109}), (\ref{eq0110}), 
(\ref{eq0301a}), (\ref{eq0301b}). Obviously, each solution of this system is a waveguide mode. 
Thus, the parameters $W, K$ should obey the dispersion equation  (\ref{eq0310}). 
Let $K$ be a variable. Equation (\ref{eq0310}) defines an explicit function $W(K)$.
Everywhere above we studied the inverse function $K(W)$ as the dispersion curve, but here the choice of $W(K)$ is more convenient. 

Fix the amplitude of the field $u(y)$, say, by setting 
\[
u(0) = 1. 
\] 
Solution $u$ depends on $y$ as on variable and on $K$ as on parameter. Thus, there 
is a family of solutions indexed by $K$. Consider the derivative 
\[
u_K (y) \equiv \ptl_K u(y). 
\]
This function obeys the equation 
\begin{equation}
u_K'' + \frac{W}{c^2} u_K - K u_K = - \frac{1}{c^2} \frac{d W}{d K} u + u
\label{eqC08}
\end{equation}
(this is (\ref{eqC01}) differentiated with respect to $K$),
and the same boundary conditions as formulated for $u$, i.~e\ (\ref{eq0109}), (\ref{eq0110}), 
(\ref{eq0301a}), (\ref{eq0301b}). 
 
Multiply (\ref{eqC08}) by $\rho u$, and integrate over the cross-section of the waveguide. 
Then multiply equation (\ref{eqC01}) for $u$ by $\rho u_K$ and integrate over the cross-section 
of the waveguide. Subtract the second integral from the first one. After taking into account the boundary conditions, obtain the relation 
\begin{equation}
\int \limits_0^{H_3} \rho(y) u^2(y) dy = 
\frac{d W}{d K}
\int \limits_0^{H_3} \frac{\rho(y)}{c^2(y)} u^2(y) dy, 
\label{eqC09}
\end{equation} 
which is valid at any point of the dispersion curve. 
Now let $W$ be a branch point of the dispersion curve, and let $u_*(y)$ be the waveguide mode
corresponding to this branch point. Due to relations (\ref{eq0203}) and (\ref{eq0204}),
$d W/ d K = 0$ at this point, thus
\begin{equation}
\int \limits_0^{H_3} \rho(y) u_*^2(y) dy = 0. 
\label{eqC10}
\end{equation} 
This is the orthogonality relation associated with the branch point. 

Note that (\ref{eqC09}) enables one to find the group velocity at any point of the 
dispersion curve as  
\begin{equation}
v_{\rm gr} =
\sqrt{\frac{K}{W}}  \int \limits_0^{H_3} \rho(y) u^2(y) dy
\left( 
\int \limits_0^{H_3} \frac{\rho(y)}{c^2(y)} u^2(y) dy
\right)^{-1}.
\label{eqC11}
\end{equation}

\section*{Appendix~B. Computation of positions of the branch points for small $\eps_1$, 
$\eps_2$}

\subsection*{B1. Classification of cases}

Here we describe how to find the approximate positions of the branch points, namely 
the values $\Theta_{\mu,\nu,m,n}(\eps_1 , \eps_2)$ and 
$\Xi_{\mu,\nu,m,n}(\eps_1 , \eps_2)$ for some small 
$\eps_1 , \eps_2$. 

Each branch points (more precisely, a pair of complex conjugate branch points) is marked 
by four indexes. They are indexes $\mu, \nu$ denoting the families of sheets and $m,n$
denoting the numbers of sheets in the families. 


{\em Case 1}
The branch point belongs to an intersection of sheets corresponding to layers 
1 and 2 or layers 2 and~3. This means that the branch point describes an interaction of synchronous waves traveling in neighboring layers. Such branch points have indexes 
$(\mu , \nu) = (1 , 2)$ or $(2 , 3)$. Note that the order of the indexes in the pair plays no role. We assume that at least one of the indexes $m,n$ is not equal to zero.

{\em Case 2}
The branch belongs to an intersection of sheets corresponding to layers 1 and 3, 
i.~e.\ $(\mu,\nu) = (1,3)$. The branch point describes an interaction between synchronous waves traveling in non-neighboring sheets. At least one index of the pair $m,n$
should be non-zero. There are two subcases that lead to different formulas:

{\em Subcase 2.1}
Both $m$ and $n$ are non-zero.

{\em Subcase 2.2}
One of the indexes $m,n$ is zero.

{\em Case 3}  
The vicinity of $K=W=0$, i.~e.\ $m = n = 0$ for any $\mu$ $\nu$. This is interaction between 
modes in all layers. This case may be of low practical importance, but formally it requires a separate consideration. In the current paper we do not build the perturbation theory for this 
case. 

Note that here we consider a general case, i.~e.\ there are no points except the origin 
where three sheets $K_{\nu ,n} (W)$ defined by (\ref{eq0306}) can intersect simultaneously.
In particular cases, when such thing happens, one should add Case~4 of three synchronous waves in all three layers. 


\subsection*{B2. Formulation of the problem in the perturbative form}

Let $\eps$ be a small parameter. Let be
\begin{equation}
\eps_1 = s_1\eps , \qquad \eps_2 = s_2 \eps,
\label{eqB01}
\end{equation} 
where $s_1$ and $s_2$ are fixed parameters. We will consider the pairs $(s_1 , s_2) = 
(1,0), (0,1), (1,1)$. The first two pairs will used for descriptions of Case~1 (linking between the neighbouring layers), the third pair will be used to describe Case~2
when all layers should be taken into account.  

Let $u(y, \eps)$ be the wave function corresponding to the branch point, i.~e.\
parameter $W$ for this solution is equal to 
$\Theta (\eps_1 , \eps_2) = \Theta_{\mu,\nu, m,n} (\eps_1, \eps_2)$, and 
parameter $K$ is equal to $\Xi (\eps_1 , \eps_2) = \Xi_{\mu,\nu,m,n} (\eps_1, \eps_2)$. 
Let the amplitude of $u(y, \eps)$ be fixed, say, by setting $u(0,\eps) =1$.

Expand all values as formal (asymptotic) power series of $\eps$: 
\begin{equation}
W = W^{(0)} + \eps W^{(1)} + (\eps)^2 W^{(2)} + \dots,
\qquad 
W^{(0)} = W_{\mu, \nu ,m,n} , 
\label{eqB02}
\end{equation}
(here $(\eps)^2$ is the second power of $\eps$), 
\begin{equation}
K = K^{(0)} + \eps K^{(1)} + (\eps)^2 K^{(2)} + \dots
\qquad 
K^{(0)} = K_{\mu, \nu,m,n},
\label{eqB03}
\end{equation}
\begin{equation}
u_* (y, \eps) = u^{(0)} (y) + \eps u^{(1)} (y)  +  \dots.   
\label{eqB04}
\end{equation}
   
Let us formulate equation and boundary conditions for the terms of the series. 
The equations  follow from (\ref{eqC01}):
\begin{equation} 
\frac{d^2 u^{(0)}}{dy^2} + \left(\frac{W^{(0)}}{c^2} - K^{(0)}\right) u^{(0)} = 0 , 
\label{eqB05}
\end{equation}
\begin{equation} 
\frac{d^2 u^{(1)}}{dy^2} + \left(\frac{W^{(0)}}{c^2} - K^{(0)}\right) u^{(1)} = 
- \left(\frac{W^{(1)}}{c^2} - K^{(1)}\right) u^{(0)} , 
\label{eqB06}
\end{equation}
Boundary conditions (\ref{eq0109}), (\ref{eq0110}) remain the same (for any asymptotic order), boundary conditions 
(\ref{eq0301a}), (\ref{eq0301b}) contain the small parameter, so they should be
 substituted by 
\begin{equation}
\frac{d u^{(1)}(H_1)}{d y} = s_1 [\rho_2 u^{(0)} (H_1 + 0) - \rho_1 u^{(0)} (H_1-0)],
\label{eqB08}
\end{equation}
\begin{equation}
\frac{d u^{(1)}(H_2)}{d y} = s_2 [\rho_3 u^{(0)} (H_2 + 0) - \rho_2 u^{(0)} (H_2-0)],
\label{eqB09}
\end{equation}


The system is insufficient to study the behavior of the branch point. It should
be completed, say, by the orthogonality relation (\ref{eqC10}), which is can be written in the form 
\begin{equation}
\int \limits_0^{H_3} \rho(y)\, (u^{(0)}(y))^2 dy = 0. 
\label{eqB12}
\end{equation}
\begin{equation}
\int \limits_0^{H_3} \rho(y)\, u^{(0)}(y) \, u^{(1)}(y) dy = 0. 
\label{eqB13}
\end{equation}


\subsection*{B3. Case 1 }

Consider the branch point $\Theta_{1,2,m,n}$. Take $s_1 = 1$, $s_2 =0$, i.~e.\ consider the
third layer to be disjoint from the first two.  
Build the zero-order approximation to the solution. The solution should have non-zero field in the first two layers and zero field in  the third layer, since the mode in the third layer is not synchronized with the first two. Since the parameters $W^{(0)}$ and $K^{(0)}$ are taken equal to 
$W_{1,2,m,n}$ and $K_{1,2,m,n}$, the field in the waveguide should be equal to 
\begin{equation}
u^{(0)} (y) = 
\left \{ \begin{array}{ll}
\cos \left( \pi m y / h_1  \right) , & 0<y<H_1 \\ 
a \cos \left( \pi n (H_2 -y) / h_2 \right) , & H_1<y<H_2 \\
0 , & H_2< y <H_3
\end{array} \right.
\label{eqB15}
\end{equation} 
for some coefficient $a$. This Ansatz guarantees that the field obeys Neumann conditions at $y = 0$
and $a = H_2$. The coefficient $a$ should be found from the orthogonality relation
(\ref{eqB12}). Namely, by substituting (\ref{eqB15}) into (\ref{eqB12}) obtain
\begin{equation}
a = a_\pm =  \pm i \sqrt{\frac{\rho_1 h_1 \sigma_1}{\rho_2 h_2 \sigma_2}}
\label{eqB16}
\end{equation}
where 
\begin{equation}
\sigma_j = \left\{ \begin{array}{ll}
2, & j = 0 \\
1, & j = 1,2,3,\dots
\end{array} \right.
\label{eqB17}
\end{equation}
There are two possible values for $a$ and they give rise to two solutions $u^{(0)}_{\pm}$. 
These solutions are zero approximations to two solutions $u_\pm (y, \eps)$, i.~e.\ there
are two solutions obeying  all conditions imposed on $u (y, \eps)$.

Our aim is to find $W^{(1)}$ and $K^{(1)}$. 
Consider the orthogonality condition (\ref{eqC04}). Substitute solutions 
$u_1 = u_+^{(0)}$, $u_2 = u_+(y, \eps)$ into this relation. Note that $u^{(0)}_+ = u_+ (y,0)$, 
thus it obeys the equation and the boundary conditions for $W_1= W^{(0)}$, 
$K_1 = K^{(0)}$. Since $\eps_1^\dag = \eps_2^\dag = 0$, 
expression (\ref{eqC07}) can be used for~$S$.    

In the order $(\eps)^0$ relation (\ref{eqC04}) is valid identically. In the order $(\eps)^1$
we get
\begin{equation}
W^{(1)} \int \limits_0^{H_3} \frac{\rho(y)}{c^2} \left( u_\pm^{(0)}(y) \right)^2 dy 
= 
\left( \rho_2 u_\pm^{(0)}(H_1 + 0) - \rho_1 u_\pm^{(0)}(H_1 - 0)\right)^2.
\label{eqB18}
\end{equation}
After substitution of (\ref{eqB15}) into (\ref{eqB18}) obtain
\begin{equation}
W^{(1)}_{\pm} = 2 \left( \sqrt{\gamma_1} \mp (-1)^{m+n} i \sqrt{\gamma_2} \right)^2
\left( 
\frac{1}{c_1^2}
-
\frac{1}{c_2^2}
\right)^{-1},
\label{eqB19}
\end{equation}
\begin{equation}
\gamma_1 = \frac{\rho_1}{\sigma_1 h_1},
\qquad
\gamma_2 = \frac{\rho_2}{\sigma_2 h_2}.
\label{eqB20}
\end{equation}
One can see that there are two positions of the branch point, they are complex conjugate, 
and the distance between the branch points of the pair is~$\sim \eps$.

Find $K^{(1)}$. For this, apply relation (\ref{eqC04}) to the solutions
$u_1 = u_\mp^{(0)}(y)$ and $u_2 = u_\pm(y , \eps)$. Consider the order of $(\eps)^1$ of the relation: 
\[
\left( 
\rho_1 u_\mp^{(0)} (H_1 - 0) - \rho_2 u_\mp^{(0)} (H_1 + 0)
\right)
\left( 
\rho_1 u_\pm^{(0)} (H_1 - 0) - \rho_2 u_\pm^{(0)} (H_1 + 0)
\right) -
\]
\begin{equation}
- W_\pm^{(1)} \int \limits_0^{H_3}
\frac{\rho(y)}{c^2(y)} u_\mp^{(0)}(y) u_\pm^{(0)}(y) dy + 
K_\pm^{(1)} \int \limits_0^{H_3}
\rho(y) u_\mp^{(0)}(y) u_\pm^{(0)}(y) dy =0. 
\label{eqB21}
\end{equation}
Since $W^{(1)}_\pm$ is already known, one can find 
\begin{equation}
K_\pm^{(1)} = 2 \frac{
\gamma_1 c_1^2 - \gamma_2 c_2^2 \mp i (-1)^{m+n} \sqrt{\gamma_1 \gamma_2} (c_1^2 + c_2^2) 
}{
c_2^2 - c_1^2
}
\label{eqB22}
\end{equation}


\subsection*{B4. Case 2}

In Case~2 the interacting layers are 1 and 3, i.~e.\
we are looking for $\Theta_{1,3,m,n}$ and~$\Xi_{1,3,m,n}$.
Since they are separated by layer 2, in which the 
the wave is not synchronous with layers 1 and 3, the link between the layers is weaker than in Case~1. Thus, the separation between the conjugate branch points is of order $(\eps)^2$, not 
~$\eps$. We need to study the second approximation of all equations to find the positions of the branch points. 

Take $s_1 = s_2 = 1$. 
Construct the zero order approximation as follows. The field in layer~2 is assumed to be zero since 
there is no phase synchronism with the other two layers. Thus, 
\begin{equation}
u^{(0)}_\pm (y) = 
\left\{ \begin{array}{ll} 
\cos (\pi m y / h_1) ,               & 0<y<H_1 \\
0 ,                                  & H_1<y<H_2 \\
a_\pm \cos (\pi n (H_3 - y) / h_3) , & H_2<y<H_3
\end{array} \right.
\label{eqB24}
\end{equation},
\[
a_\pm =\pm i \sqrt{\frac{h_1 \rho_1 \sigma_1}{h_3 \rho_3 \sigma_3}}
\]
 
Take (\ref{eqC04}) with $u_1 = u^{(0)}_\pm (y)$ and $u_2 = u_\pm (y, \eps)$. 
In the first approximation (similarly to the previous subsection) obtain 
\begin{equation}
W^{(1)}_\pm = W^{(1)} = 2 (\gamma_1 - \gamma_3) \left( \frac{1}{c_1^2} - 
\frac{1}{c_3^2} \right)^{-1},
\label{eqB25}
\end{equation} 
\[
\gamma_1 = \frac{\rho_1}{\sigma_m h_1},
\qquad
\gamma_3 = \frac{\rho_3}{\sigma_n h_3}.
\]
Take (\ref{eqC04}) with $u_1 = u^{(0)}_\pm (y)$ and $u_2 = u_\mp (y, \eps)$. 
In the first approximation get
\begin{equation}
K^{(1)}_\pm = K^{(1)} = 2 \frac{\gamma_1 c_1^2 - \gamma_3 c_3^2}{c_3^2 - c_1^2}. 
\label{eqB26}
\end{equation} 
One can see that in the first approximation the branch points are not separated. 
So we should construct $u^{(1)}_\pm$ (we do it in the next subsections) and study the second approximation of (\ref{eqC04}). 
Let $u_\pm^{(1)} (y)$ be found. Substituting 
$u_1 = u^{(0)}_\pm (y)$ and $u_2 = u_\pm (y, \eps)$ into (\ref{eqC04})
and taking the terms $~\eps^2$, obtain 
\[
W^{(2)}_{\pm} = \frac{2}{\rho_1 h_1} \left( \frac{1}{c_1^2}-\frac{1}{c_3^2} \right)^{-1}
\left[
-W^{(1)} \int \limits_0^{H_3} 
\frac{\rho(y)}{c^2(y)}\, u^{(0)}_\pm (y)\, u^{(1)}_\pm (y) \, dy
 + \right.
\]
\[
(-1)^n \rho_3 a_\pm (\rho_3 u^{(1)}_\pm (H_2 + 0)-\rho_2 u^{(1)}_\pm (H_2 - 0)) 
- 
\qquad \qquad \qquad \qquad \qquad  
\].
\begin{equation} 
\qquad \qquad \qquad \qquad \qquad 
\left.
 (-1)^m \rho_1 (\rho_2 u^{(1)}_\pm (H_1 + 0)-\rho_1 u^{(1)}_\pm (H_1 - 0)) 
\right] 
\label{eqB27}
\end{equation}

Then, substituting 
$u_1 = u^{(0)}_\mp (y)$ and $u_2 = u_\pm (y, \eps)$ into (\ref{eqC04}),
and combining the result with (\ref{eqB27}), obtain
an expression for $K^{(2)}_{\pm}$.
After some algebra obtain
\[
K_\pm^{(2)} = \frac{W^{(2)}}{c_1^2} + 
\frac{4 \gamma_1}{h_1 \sigma_1} \int \limits_0^{H_1} u^{(0)}_\pm (y) u^{(1)}_\pm (y)\, dy + 
\]
\begin{equation}
+\frac{2 (-1)^m}{h_1 \sigma_1} [\rho_2 u^{(1)}_\pm (H_1+0) - \rho_1 u^{(1)}_\pm (H_1 -0)].
\label{eqB27a}
\end{equation}

{\bf  Subcase 2.1}

Let be $m \ne 0$, $n \ne 0$. 
Construct solution $u^{(1)}_\pm$ obeying equation (\ref{eqB06}), boundary conditions 
(\ref{eq0109}), (\ref{eq0110}), (\ref{eqB08}), (\ref{eqB09}), and relation 
$u^{(1)}_\pm (0) =0$ for uniqueness.
This problem has the following solution:

For $0<y<H_1$:
\begin{equation}
u^{(1)}_\pm (y) = - \frac{\rho_1 y}{\pi m} \sin \left( 
\frac{\pi m}{h_1} y
\right).
\label{eqB29}
\end{equation} 
For $H_1 < y < H_2$: 
\begin{equation}
u^{(1)}_\pm (y) =
- \frac{(-1)^m \rho_1 
\cos(\alpha_2^{(0)} (y-H_2))
+
(-1)^n a_\pm \rho_3
\cos(\alpha_2^{(0)} (y-H_1))
}{\alpha_2^{(0)} \sin(\alpha_2^{(0)} h_2)} ,
\label{eqB30}
\end{equation}
where 
\begin{equation}
\alpha_2^{(0)} = \sqrt{\frac{W^{1,3}_{m,n}}{c_2^2} - K^{1,3}_{m,n}} .
\label{eqB30}
\end{equation}
For $H_2 < y < H_3$: 
\begin{equation}
u^{(1)}_\pm (y) = - 
\frac{a_{\pm} \rho_3 (H_3 -y)}{\pi n} \sin \left( 
\frac{\pi n}{h_3} (H_3-y)
\right)
+
\beta a_\pm
\cos \left( 
\frac{\pi n}{h_3} (H_3-y)
\right),
\label{eqB29a}
\end{equation} 
coefficient $\beta$ is found from relation (\ref{eqB13}):
\begin{equation}
\beta = \frac{\rho_1 h_1}{2 \pi^2 m^2} - \frac{\rho_3 h_3}{2 \pi^2 n^2}. 
\label{eqB31}
\end{equation}

Substitute $u^{(1)}_\pm$ into 
(\ref{eqB27}) and (\ref{eqB27a}).
The result is as follows: 
\[
W^{(2)}_\pm = 
\left( \frac{1}{c_1^2} - \frac{1}{c_3^2} \right)^{-1}
\left[
\pi^{-2}
\left( 
\frac{\rho_3^2}{n^2} 
-
\frac{\rho_1^2}{m^2}
\right)
+ 
\frac{2 \rho_2 \cos(\alpha_2^{(0)} h_2) (\gamma_1 -\gamma_3)}{
\alpha_2^{(0)} \sin (\alpha_2^{(0)} h_2)}
\right. 
\]
\begin{equation}
\left.
\qquad \qquad \qquad \qquad \qquad \qquad 
\pm\frac{4 i (-1)^{m+n} \rho_2 \sqrt{\gamma_1 \gamma_3}}{
\alpha_2^{(0)} \sin (\alpha_2^{(0)} h_2)
}
\right],
\label{eqB32}
\end{equation} 

\[
K^{(2)}_\pm  = 
\frac{1}{c_1^2 -c_3^2} \left[ 
\pi^{-2} \left( \frac{c_1^2 \rho_1^2}{m^2}  - \frac{c_3^2 \rho_3^2}{n^2} \right) 
+ \frac{2 \rho_2 (\gamma_3 c_3^2 - \gamma_1 c_1^2) \cos(\alpha_2^{(0)} h_2)}{
\alpha_2^{(0)} \sin (\alpha_2^{(0)} h_2)
}
\right.
\]
\begin{equation}
\left.
\qquad \qquad \qquad \qquad \qquad 
\mp 
\frac{
2 i (-1)^{m+n} \rho_2 (c_1^3 + c_3^2)  \sqrt{\gamma_1 \gamma_3} 
}{
\alpha_2^{(0)} \sin (\alpha_2^{(0)} h_2)
}
\right].
\label{eqB33}
\end{equation}


{\bf Subcase 2.2}

Let be $m = 0$, $n \ne 0$. 
The first-order approximation of the solution has the following form.
For $0 < y < H_1$: 
\begin{equation}
u^{(1)}_\pm (y) = - \frac{\rho_1 y^2}{2 h_1},
\label{eqB34}
\end{equation}
for $H_1 < y < H_2$: 
\begin{equation}
u^{(1)}_\pm (y) = -
\frac{
2 \rho_1 \cos( \alpha_2^{(0)} (y-H_2)) + (-1)^n a_{\pm} \rho_3 \cos(\alpha_2^{(0)} (y - H_1))
}{
\alpha_2^{(0)} \sin(\alpha_2^{(0)} h_2), 
}
\label{eqB35}
\end{equation}
For $H_2 < y < H_3$: 
\begin{equation}
u^{(1)}_\pm (y) = 
- \frac{a_\pm \rho_3 (H_3 - y)}{\pi n} \sin \left( 
\frac{\pi n}{h_3} (H_3 - y)
\right)
+ 
\beta a_{\pm} \cos \left( 
\frac{\pi n}{h_3} (H_3 - y)
\right),
\label{eqB36}
\end{equation}
\begin{equation}
\beta = -\frac{h_1 \rho_1}{6} - \frac{\rho_3 h_3}{2 \pi^2 n^2}.
\label{eqB37}
\end{equation}

Using this solution similarly to Subcase~2.1, obtain
\[
W^{(2)}_\pm = 
\left( \frac{1}{c_1^2} - \frac{1}{c_3^2} \right)^{-1}
\left[
\left( 
\frac{\rho_3^2}{\pi^2 n^2} 
-
\frac{\rho_1^2}{3}
\right)
+ 
\frac{2 \rho_2 \cos(\alpha_2^{(0)} h_2) (\gamma_1 -\gamma_3)}{
\alpha_2^{(0)} \sin (\alpha_2^{(0)} h_2)}
\right. 
\]
\begin{equation}
\left.
\qquad \qquad \qquad \qquad \qquad \qquad 
\pm\frac{4 i (-1)^{m+n} \rho_2 \sqrt{\gamma_1 \gamma_3}}{
\alpha_2^{(0)} \sin (\alpha_2^{(0)} h_2)
}
\right],
\label{eqB38}
\end{equation} 

\[
K^{(2)}_\pm  = 
\frac{1}{c_1^2 -c_3^2} \left[ 
\left( \frac{c_1^2 \rho_1^2}{3}  - \frac{c_3^2 \rho_3^2}{\pi^2 n^2} \right) 
+ \frac{2 \rho_2 (\gamma_3 c_3^2 - \gamma_1 c_1^2) \cos(\alpha_2^{(0)} h_2)}{
\alpha_2^{(0)} \sin (\alpha_2^{(0)} h_2)
}
\right.
\]
\begin{equation}
\left.
\qquad \qquad \qquad \qquad \qquad 
\mp 
\frac{
2 i (-1)^{n} \rho_2 (c_1^3 + c_3^2)  \sqrt{\gamma_1 \gamma_3} 
}{
\alpha_2^{(0)} \sin (\alpha_2^{(0)} h_2)
}
\right].
\label{eqB39}
\end{equation}


\end{document}